\documentclass{amsart} 
\usepackage[normalem]{ulem}
\usepackage{amsmath,amsthm,amsfonts,amscd,amssymb, tikz-cd} 
\usepackage{verbatim}      	
\usepackage{listings} 
\usepackage{mathtools}
\usepackage{epsfig}         	
\usepackage{url}		
\usepackage{epic,eepic,latexsym}
\usepackage{graphicx}
\usepackage{xcolor}

\usepackage{lcd}
\usepackage{mathrsfs}
\usepackage[centering,text={15.5cm,22cm},
        marginparwidth=20mm,dvips]{geometry}
\usepackage[linktocpage]{hyperref}
\usepackage{lscape}
\usepackage{tabularx}
\usepackage{marginnote}
\usepackage[all]{xy}
\usepackage{enumitem}
\usepackage{stmaryrd}
\usepackage{tipa}
\usepackage{upgreek}

\usetikzlibrary {arrows.meta,positioning}
\usepackage{caption}
\usepackage{subcaption}

\DeclareMathOperator{\tw}{tw}

\newcommand{\QQ}{\mathbb{Q}}

\DeclareMathOperator{\sign}{sgn}

\DeclareMathOperator{\prin}{prin}

\DeclareMathOperator{\Li}{Li}

\DeclareMathOperator{\uf}{uf}

\DeclareMathOperator{\Sk}{Sk}

\def\!{\mskip-\thinmuskip}
\let\?=\overline
\let\:=\colon
\let\llb=\llbracket
\let\rrb=\rrbracket
\let\bb=\mathbb

\let\rar=\rightarrow
\let\f=\mathfrak

\let\wh=\widehat
\let\wt=\widetilde

\newcommand {\kk} {\Bbbk}

\newcommand {\mm} {{\bf m}}

\newcommand {\SSS} {{\bf S}}
\newcommand {\sd} {{\bf s}}

\theoremstyle{plain}
\ifdefined\thm
\else
	\theoremstyle{plain}
 	\newtheorem{thm}{Theorem}[section]
\fi

\ifdefined\prop
\else
	\theoremstyle{plain}
 	\newtheorem{prop}[thm]{Proposition}
\fi\ifdefined\cor
\else
	\theoremstyle{plain}
 	\newtheorem{cor}[thm]{Corollary}
\fi

\ifdefined\lem
\else
	\theoremstyle{plain}
 	\newtheorem{lem}[thm]{Lemma}
\fi

\ifdefined\def
\else
	\theoremstyle{definition}
	\newtheorem{def}[thm]{Definition}
\fi

\ifdefined\defn
\else
	\theoremstyle{definition}
	\newtheorem{defn}[thm]{Definition}
\fi

\ifdefined\example
\else
	\theoremstyle{definition}
    \newtheorem{example}[thm]{Example}
\fi

\ifdefined\eg
\else
		\theoremstyle{definition}
        \newtheorem{eg}[thm]{Example}
\fi

\ifdefined\rem
\else
		\theoremstyle{remark}
        \newtheorem{rem}[thm]{Remark}
\fi

\ifdefined\quest
\else
		\theoremstyle{plain}
		\newtheorem{quest}[thm]{Question}
\fi

\ifdefined\assumption*
\else
		\theoremstyle{plain}
		\newtheorem*{assumption*}{Assumption}
\fi

\newcommand{\Q}{{\mathbb Q}}
\newcommand{\R}{{\mathbb R}}
\newcommand{\Z}{{\mathbb Z}}

\newcommand{\fD}{\mathfrak{D}}

\newcommand{\Gr}{\operatorname{Gr}}

\newcommand{\supp}{\operatorname{supp}}

\usepackage{tikz-cd}

\newcommand{\incoming}{\mathrm{in}}
\newcommand{\Quot}{\mathrm{Quot}}

\newcommand{\frg}{\mathfrak{g}}
\newcommand{\frD}{\mathfrak{D}}

\newcommand{\frp}{\mathfrak{p}}
\newcommand{\frd}{\mathfrak{d}}

\newcommand{\hfrg}{\hat{\frg}}
\newcommand{\hG}{\hat{G}}

\newcommand{\cM}{M^{\circ}}
\newcommand{\bJ}{\overline{J}}

\newcommand{\bu}{\overline{u}}
\newcommand{\bv}{\overline{v}}
\newcommand{\bV}{\overline{V}}
\newcommand{\bU}{\overline{U}}
\newcommand{\bA}{\overline{A}}
\newcommand{\bQ}{\overline{Q}}
\newcommand{\bd}{\overline{d}}
\newcommand{\bW}{\overline{W}}
\newcommand{\bR}{\overline{R}}
\newcommand{\bsig}{\overline{\sigma}}

\newcommand{\bB}{\overline{B}}

\newcommand{\bI}{\overline{I}}
\newcommand{\bP}{\overline{P}}

\newcommand{\bn}{\overline{n}}
\newcommand{\bi}{\overline{i}}
\newcommand{\bj}{\overline{j}}
\newcommand{\bk}{\overline{k}}
\newcommand{\bp}{\overline{p}}
\newcommand{\ba}{\overline{a}}

\newcommand{\bM}{\overline{M}}
\newcommand{\bN}{\overline{N}}
\newcommand{\bcM}{\overline{M}^\circ}

\newcommand{\st}{\mathrm{st}}

\newcommand{\tQ}{\tilde{Q}}

\newcommand{\NN}{\mathbb{N}}
\newcommand{\ZZ}{\mathbb{Z}}

\newcommand{\op}{{\mathrm{op}}}

\newcommand{\ufv}{{\operatorname{uf}}}
\newcommand{\fv}{{\mathrm f}}

\renewcommand{\tw}{\operatorname{tw}}

\ifdefined\supp
\else
\newcommand{\supp}{\operatorname{supp}}
\fi

\ifdefined\tw
\else
	\newcommand{\tw}{{\opname tw}}
\fi

\ifdefined\G
    \renewcommand{\G}{{\mathbb G}}
\else
   \newcommand{\G}{{\mathbb G}}
\fi

\ifdefined\C
    \renewcommand{\C}{{\mathbb C}}
\else
   \newcommand{\C}{{\mathbb C}}
\fi

\usepackage{comment}

\author{Qiyue Chen}
\address{School of Mathematics \\ University of Minnesota \\ Minneapolis, MN 55455\\ USA}
\email{chen8448@umn.edu }

\author{Travis Mandel}
\address{Department of Mathematics \\ University of Oklahoma \\ Norman, OK 73019 \\ USA}
\email{tmandel{\char'100}ou.edu}

\author{Fan Qin}
\address{School of Mathematical Sciences \\ Beijing Normal University \\ China}
\email{qin.fan.math@gmail.com}

\title{Stability scattering diagrams and quiver coverings}

\begin{document}
\begin{abstract}
Given a covering of a quiver (with potential), we show that the associated Bridgeland stability scattering diagrams are related by a restriction operation under the assumption of admitting a nice grading.  We apply this to quivers with potential associated to marked surfaces.  In combination with recent results of the second and third authors, our findings imply that the bracelets basis for a once-punctured closed surface coincides with the theta basis for the associated stability scattering diagram, and these stability scattering diagrams agree with the corresponding cluster scattering diagrams of Gross-Hacking-Keel-Kontsevich except in the case of the once-punctured torus.
\end{abstract}


\maketitle
\setcounter{tocdepth}{1} \tableofcontents{}
	
\section{Introduction}	\label{sec:intro}

\subsection*{Background}

Let $S$ denote a compact oriented surface and $M$ a finite set of marked points in $S$. The pair $\Sigma=(S, M)$ is called a marked surface. Marked points that are located within the interior of $S$ are known as punctures. We can assume $S$ is connected without loss of generality. In the fields of cluster algebras, tau-tilting theory, and related areas, it is more challenging to study once-punctured closed surfaces $\Sigma$ than it is to study other types of surfaces.

For example, by \cite[Thm. 1.2.4]{qin2019bases} or \cite[Cor. 1.2]{mou2019scattering}, the cluster scattering diagrams in \cite{gross2018canonical} and Bridgeland's stability scattering diagrams \cite{bridgeland2017scattering} are equivalent in many cases. The former diagram has a powerful combinatorial construction, while the latter encodes important information about semistable objects appearing in the corresponding representation theory. Notably, 
in combination with \cite[Thm. 1.2]{yurikusa2020density}, these results imply this equivalence of scattering diagrams for all marked surfaces expect the once-punctured closed surfaces. It is natural to inquire about a comparison for these exceptional cases.  
\begin{quest}\label{quest:diagrams}
Are the cluster scattering diagram and the stability scattering diagram equivalent for once-punctured closed surfaces?
\end{quest}

Furthermore, by studying the topology (or Teichmüller theory) of $\Sigma$, one can construct the (tagged) skein algebra $\Sk(\Sigma)$. This algebra has a basis consisting of the (tagged) bracelets, which can be constructed as topological diagrams \cite{frohman2000skein,fock2006moduli, musiker2013bases}. By \cite{MandelQin2021}, if $\Sigma$ is NOT a once-punctured torus, the bracelets coincide with the theta functions of the related cluster scattering diagrams in the sense of \cite{gross2018canonical}. These theta functions are crucial constructions in recent advancements of cluster theory.
\begin{quest}\label{quest:theta-tori}
    Can we interpret the (tagged) bracelets for once-punctured tori as theta functions for the stability scattering diagram?
\end{quest}

\subsection*{Main Result}

Let $\pi: Q\rightarrow \bQ$ be a covering map between quivers, see \S \ref{sec:covering_representation} (this roughly means that $\pi$ is the projection modulo the action of a finite group $\Pi$). These are the skew-symmetric cases of the coverings considered in  \cite{felikson2012cluster,huang2018unfolding}. We have an algebra homomorphism $\sigma:\widehat{\C \bQ}\rightarrow \widehat{\C Q}$ between the completed path algebras which sends a path $p$ to the sum of the paths in $\pi^{-1}(p)$.  Following \cite{DerksenWeymanZelevinsky08}, given potentials $W$ and $\bW$ for $Q$ and $\bQ$, respectively, we mod out by the associated Jacobian ideals to construct the completed Jacobian algebras $J$ and $\bJ$, respectively.  We assume $W=\sigma(\bW)$.

Now $\sigma$ descends to an algebra homomorphism $\bsig:\bJ\rar J$, and any $J$-module $V$ pulls back to a $\bJ$-module $\bsig^* V$.  Let $P_k$ denote the $k$-th indecomposable projective left module $J e_k$ of $J$ and $\bP_{\bk}$ the $\bk$-th indecomposable projective left module $\bJ e_{\bk}$ of $\bJ$, respectively, where $\bk\coloneqq\pi(k)$.  Let $\fD^\st(J)$ and $\fD^\st(\bJ)$ denote the stability scattering diagrams \cite{bridgeland2017scattering} associated to $J$ and $\bJ$, respectively. We use the motivic Hall algebras of nilpotent\footnote{We emphasize the nilpotency condition since it is needed for Lemma \ref{lem:stability_projective} when the Jacobian algebra is infinite-dimensional.  However, as pointed out to us by Bernhard Keller, all finite-dimensional representations of a \textit{completed} Jacobian algebra are automatically nilpotent; cf. \cite[\S 10; argument due to Bill Crawley-Boevey]{DerksenWeymanZelevinsky08}.  So the nilpotency condition is redundant in this paper.} representations when constructing these diagrams, see \cite[\S 7]{Nagao10}. The \textit{restriction} $\?{\fD^\st(J)}$ is a scattering diagram roughly obtained by taking a ``slice'' of the walls of $\fD^\st(J)$ (i.e., intersecting their supports with a certain subspace) and applying a certain projection map to the scattering functions of $\fD^\st(J)$; see \S \ref{sec:restriction} for details.

\begin{thm}\label{thm:stability_diagram}
    Let $\pi:(Q,W)\rightarrow (\bQ,\bW)$ be a covering as above.  If every indecomposable projective module of $J$ has a nice grading, then the restriction $\?{\fD^\st(J)}$ is equivalent to $\fD^\st(\bJ)$.
\end{thm}

The nice grading in Theorem \ref{thm:stability_diagram} is introduced in \S \ref{sec:covering_representation} as a special instance of the grading in \cite{haupt2012euler}, but we will construct nice grading for modules $P_k$ which are often not the tree modules or band modules studied in  \cite{haupt2012euler}. Briefly, a nice grading for $P_k=Je_k$ associates an integer $\partial(V)$ to each vertex $V$ in the support of $P_k$, inducing $\partial(a)\coloneqq \partial(t(a))-\partial(s(a))$ for arrows $a$, subject to two conditions: (1) distinct vertices in the same fiber of $\pi$ must map to \textit{different} integers, and (2) two arrows in the support of $P_k$ which lie in the same fiber of $\pi$ must map to the \textit{same} integer. 
 We give more intuitive sufficient conditions below in Theorem \ref{thm:intro-non-wrap} and Corollary \ref{cor:intro-cover}. See \cite{zhou2020cluster,MandelQin2021} for related results on restrictions of cluster scattering diagrams under quiver folding.\footnote{As in \cite[\S 4]{felikson2012cluster} and \cite[\S 2]{huang2018unfolding}, a covering is called a folding if it satisfies some additional technical conditions.}

We apply Theorem \ref{thm:stability_diagram} to surfaces where the covering map arises from topological constructions.  Let $\Sigma$ be a connected once-punctured closed surface $(S,M)$, and let $\sd$ be a seed associated to an ideal triangulation $T$ of $\Sigma$ without frozen vertices. Let $Q$ denote the quiver associated to $T$, $Q^{\op}$ the opposite quiver,\footnote{The opposite quiver $Q^{\op}$ appears by our choice of conventions. We do not need to take the opposite if we attach a negative sign to the matrix $B$ associated to the the quiver $Q$ in \S \ref{sec:diagrams}.} and $J^{\op}$ the completed Jacobian algebra associated to $Q^{\op}$ with potential given by \cite[Definition 23]{Lad09}. Let $\fD(\sd)$ denote the corresponding cluster scattering diagram \cite{gross2018canonical}, see \S \ref{sec:diagrams}.

\begin{thm}\label{thm:1p_surface_diagrams}
For $J^{\op}$ and $\sd$ associated to a once-punctured surface $(S,M)$ as above, $\fD^\st(J^{\op})$ and $\fD(\sd)$ are equivalent if and only if $S$ is not a torus.
\end{thm}
Indeed, the equivalence extends to arbitrary coefficients, see Remark \ref{rem:equivalence_coefficinets}.

The equivalence between the two scattering diagrams for non-degenerate potentials was previously known under the injective-reachability condition \cite{qin2019bases,mou2019scattering} (equivalently, when there exists a green-to-red sequence \cite{keller2011cluster,muller2015existence}). Consequently, the equivalence was known for all triangulable marked surfaces except for the once-punctured closed surfaces.  Theorem \ref{thm:1p_surface_diagrams} addresses these exceptional cases and answers Question \ref{quest:diagrams}. Notably, it extends the equivalence to cases when the injective-reachability condition is not satisfied. In addition, we obtain an intriguing example (the torus case) when the two scattering diagrams are NOT equivalent.\footnote{It was shown in \cite[Cor. 1.2(ii)]{mou2019scattering} that the cluster and stability scattering diagrams in the once-punctured torus case differ along at most a single wall.  Our result confirms that this difference is nontrivial.}

In \cite{MandelQin2021}, it was proved that the tagged bracelets basis for the tagged skein algebra coincides with the theta basis constructed from $\fD(\sd)$, except for the once-punctured torus, where one should instead use the restriction $\?{\fD(\wt{\sd})}$ for $\wt{\sd}$ the seed associated to a covering space of the punctured torus (i.e., a torus with multiple punctures). Our proof of Theorem \ref{thm:1p_surface_diagrams} shows that $\?{\fD(\wt{\sd})}$ is equivalent to $\fD^\st(J^{\op})$; cf. Lemma \ref{lem:cover-stability}. Thus, we find that the tagged bracelets basis always agrees with the theta basis constructed from the stability scattering diagram $\fD^\st(J^{\op})$; i.e., we have an affirmative answer for Question \ref{quest:theta-tori}. 

\begin{thm}\label{thm:theta-stability}
    For any triangulable marked surface $\Sigma$ with the potential given by \cite{Lad09}, the tagged bracelets basis coincides with the theta basis for the stability scattering diagram.
\end{thm}

Theorem \ref{thm:theta-stability} also has philosophical importance.  It suggests that the stability scattering diagram, rather than the cluster scattering diagram, is the correct scattering diagram to use for the construction of some canonical bases of interest.

The proofs of Theorems \ref{thm:1p_surface_diagrams} and \ref{thm:theta-stability} are based on Lemma \ref{lem:surface_nice_grading}, which says that for $d:1$ cyclic coverings of marked surfaces ($d\geq 3$), the associated covering of quivers with potential admits a nice grading.  In \S \ref{sec:gen}, we abstractify and generalize the proof of Lemma \ref{lem:surface_nice_grading} to obtain Theorem \ref{thm:non-wrap}.  For coverings with cyclic deck group (or solvable deck group; cf. Remark \ref{rem:solvable}), this result replaces the technical ``nice grading'' condition of Theorem \ref{thm:stability_diagram} with the more intuitive condition that $W$ is ``non-wrapping.''

\begin{thm}[Theorem \ref{thm:non-wrap}]\label{thm:intro-non-wrap}
    Let $\pi:(Q,W)\rar (\bQ,\bW)$ be a covering of quivers with potential with deck group $\Pi=\Z/d\Z$.  If $W$ is non-wrapping, then $\?{\f{D}^{\st}(J)}$ is equivalent to $\f{D}^{\st}(\bJ)$.
\end{thm}

The non-wrapping condition is explained as follows.  Say we partition $Q$ into ``sheets'' for the covering map $\pi$, labelled by elements of $\Z/d\Z$ which are identified with vertices of a $d$-gon in the natural way.  Each arrow in $W$ can be viewed as going up or down some number of sheets, i.e., wrapping clockwise or counterclockwise around this $d$-gon.  For each arrow, we make a choice between clockwise and counterclockwise.  Non-wrapping means that these choices can be made in a way that gives each term of $W$ a winding number of $0$ when viewed as wrapping around this polygon.

We then apply Theorem \ref{thm:intro-non-wrap} in Examples \ref{eg:general-manifold-covering} and \ref{eg:Lie-Grass}, the former of which is summarized as the following significant generalization of Lemma \ref{lem:surface_nice_grading}.

\begin{cor}\label{cor:intro-cover}
    Consider $(\bQ,\bW)$ a quiver with potential.  Suppose the topological realization of $\bQ$ is immersed in a manifold with boundary $\?{\SSS}$ such that, for each term of $\bW$, the associated loop in $\?{\SSS}$ is contractible.  Then for any normal covering space $\pi:\SSS\rar \?{\SSS}$ with solvable deck group, the induced covering $\pi:(Q,W)\rar (\bQ,\bW)$ of quivers with potential is such that $\?{\f{D}^{\st}(J)}$ is equivalent to $\f{D}^{\st}(\bJ)$.
\end{cor}

\subsection*{Acknowledgments}
The authors thank Min Huang and Bernhard Keller for inspiring discussions. The third author was supported by the National Natural Science Foundation of China (Grant No. 12422102).

\section{Scattering Diagrams}\label{sec:scattering_diagram}

Let us recall the definition of scattering diagrams \cite{gross2018canonical} following the convention of \cite{qin2019bases}. We shall use the subscript $_\R$ to denote an extension to $\R$. Let $\langle\ ,\ \rangle$ denote the natural pairing between dual spaces.

\subsection*{Seeds}\label{sec:seeds}

Let $I$ denote a finite set of vertices with a chosen partition $I=I_{\ufv}\sqcup I_{\fv}$ dividing it into unfrozen and frozen vertices. Let $d_i$ denote strictly positive rational numbers for each $i\in I$.

Let $N$ denote a rank-$|I|$ lattice with a chosen basis $E=(e_i)_{i\in I}$. Let $M$ denote the dual lattice with the dual basis $(e_i^*)_{i\in I}$. We further denote $F=(f_i)_{i\in I}$ such that $f_i=\frac{1}{d_i}e^*_i\in M_\R$. Let $\{\ ,\ \}$ denote a $\Q$-valued skew-symmetric bilinear form on $N$. We define the bilinear form $B$ on $N$ such that
\begin{align*}
B(e_i,e_j)\coloneqq B_{ij}\coloneqq \{e_j,e_i\}d_i.
\end{align*}
We further require that $B_{ij}\in \Z$ whenever $i$ or $j$ is unfrozen.

We use $\sd$ to denote the collection $(I_{\ufv},I_{\fv},(d_i)_{i\in I}, N,E,M,B)$ and call it a seed. Denote $N_{\ufv}=\bigoplus_{k\in I_\ufv} \Z e_k$ and $\cM=\bigoplus_{i\in I} \Z f_i\subset M_{\QQ}$. Define the linear map $p^*:N_{\ufv}\rightarrow \cM$ by $p^*(n)=\{n,\ \}$, i.e., \begin{align*}
p^*(e_k)=\sum_{i\in I} B_{ik}f_i.
\end{align*}

We say $\sd$ satisfies the injectivity assumption if $p^*$ is injective. When this assumption fails, we construct a new seed $\sd^{\prin}$ satisfying this assumption. More precisely, for each vertex $i\in I$, we append a new vertex $i'$ with $d_{i'}=1$. We define $\{e_i,e_{j'}\}=\delta_{ij}$, $\{e_{i'},e_{j'}\}=0$, and the rest of the data is extended in the obvious way. $\sd^{\prin}$ is called the principal coefficient seed associated to $\sd$, see \cite[\S 2.3]{MandelQin2021} for details.

Notice that $B$ is skew-symmetric if and only if all $d_i$ are equal. In this case, we associate to $B$ a quiver $\tQ$ without loops or oriented $2$-cycles such that 
\begin{itemize}
\item Its vertex set $\tQ_0$ is given by $I$,
\item $B_{ij}$ equals the number of arrows from $i$ to $j$ minus the number of arrows from $j$ to $i$.
\end{itemize}
Note that $\tQ$ comes equipped with a partition of its vertices as either frozen or unfrozen.  For our purposes, all quivers will be assumed to have such a partition of their vertices.

Conversely, any such quiver determines a seed $\sd$.  Such a seed is said to be skew-symmetric or of quiver type. Let $\tQ_1$ denote the set of arrows of $\tQ$, and $\tQ_0=I$ the set of vertices.

We use $Q$ denote the full subquiver of $\tQ$ whose vertices are unfrozen.

\subsection*{Lie algebras}\label{sec:Lie_alg}
Fix a field $\kk$ of characteristic $0$. Define the group ring $\kk[N_{\ufv}]=\bigoplus_{n\in N_{\ufv}} \kk y^n$ and 
$\kk[\cM]=\bigoplus_{m\in \cM} \kk x^m$. Denote $y_k\coloneqq y^{e_k}$ and $x_i=x^{f_i}$.

Denote the submonoid $N^\oplus\coloneqq \{\sum_{k\in I_{\ufv}} n_k e_k|n_k\geq 0\}$ and $N^+=N^\oplus \setminus\{0\}$. Define $M^\oplus=p^*(N^\oplus)$  and $M^+=p^*(N^+)$. Let $\kk\llb N^\oplus \rrb$ denote the completion of $\kk[N^\oplus]$ with respect to its maximal ideal $\kk[N^+]$ and similarly for the completion $\kk \llb M^\oplus \rrb$ of $\kk[M^{\circ}]$ with respect to $\kk[M^+]$. We introduce the formal completions: 
\begin{align*}
\kk\llb N_{\uf} \rrb&=\kk[N_{\uf}]\otimes_{\kk[N^\oplus]} \kk\llb N^\oplus \rrb,\\
\kk\llb M^{\circ} \rrb&=\kk[M^{\circ}]\otimes_{\kk[M^\oplus]} \kk\llb M^\oplus \rrb.
\end{align*}

We endow $\kk\llb N_{\uf}\rrb$ with the Poisson bracket such that $\{ y^n,y^{n'}\}=-\{n,n'\} y^{n+n'}$. Then it becomes an $N$-graded\footnote{Here and elsewhere, when we say ``graded'' or ``spanned,'' we really mean \textit{topologically} graded or spanned, essentially meaning that we allow for linear combinations to be infinite (more precisely, we mean the underlying module is the closure of a usual graded module or span with respect to the relevant adic topology).    See \cite[\S 2.2.2]{davison2019strong} for details.} Poisson algebra. Let $\frg$ denote its $N^+$-graded Lie subalgebra. 

Introduce a $\Z$-valued function $|\ |$ on $N_{\ufv}$ such that $|\sum n_k e_k|=\sum n_k$. For any $l\in \NN$, let $\frg_{>l}$ denote the Lie algebra ideal spanned by the monomials $y^n$ with $|n|>l$ and define the nilpotent Lie algebra $\frg_l=\frg/\frg_{>l}$. Let $G_l$ denote the pro-unipotent group $\exp \frg_l$ defined via the Baker-Campbell-Hausdorff formula. Let $\hfrg$ denote the inverse limit of $\frg_l$ and $\hG$ the inverse limit of $G_l$. We extend the bijection $\exp: \frg_l\simeq G_l$ to $\exp:\hfrg\simeq \hG$. Denote the natural projections $\pi_{l}:\hfrg\rightarrow \frg_l$ and $\pi_{l}:\hG\rightarrow G_l$.

For any element $n\in N^+$, we define $\hfrg_{n}^\parallel\coloneqq \kk \llb y^n\rrb$ which is an abelian Lie subalgebra of $\hfrg$. Let $\hG_n^\parallel=\exp (\hfrg_{n}^\parallel)$ denote the corresponding pro-unipotent group. It is an abelian subgroup of $\hG$.

Finally, the Poisson algebra $\hfrg$ linearly acts on $\kk\llb \cM \rrb$ via the derivations $\{\kk[N^+],\ \}$ defined by
\begin{align*}
\{y^n,x^m\}\coloneqq \langle n,m\rangle x^{m+p^*(n)}.
\end{align*}
This induces an action of $\hG$ on $\kk\llb \cM \rrb$ which is faithful by the injectivity of $p^*$.

\subsection*{Scattering Diagrams}\label{sec:diagrams}
A wall is a pair $(\frd,\frp_{\frd})$ such that:
\begin{itemize}
\item $\frd$ is a codimension-$1$ polyhedral cone in $M_\R$ contained in $(n_0)^\perp$ for some primitive vector $n_0\in N^+$.
\item $\frp_{\frd}\in \hG^\parallel_{n_0}$.
\end{itemize}
We call $\frd$ its support and $\frp_{\frd}$ its wall-crossing operator. The wall is called trivial if $\frp_{\frd}$ is the identity. It is said to be incoming if $p^*(n_0)\in \frd$ and outgoing otherwise. When the context is clear, we denote the wall simply by $\frd$.

Let $\frD$ denote a collection of walls $(\frd,\frp_{\frd})$. We use $\frD_l$ to denote the subset of walls $(\frd,\frp_{\frd})\in \frD$ with $\pi_l(\frp_{\frd})\neq 0$. We say $\frD$ is a scattering diagram if $\frD_l$ only consists of finitely many walls for each $l\in \NN$. Its (essential) support $\supp \frD$ is defined as the union of the supports of its non-trivial walls. The closures of the $|I|$-dimensional connected components of $M_{\bb{R}}\setminus \?{\supp \frD}$ are called chambers. The joints of $\frD$ are the relative boundary components of walls and the codimension-$2$ intersections of walls $\frd_1\cap \frd_2$.

A smooth path $\gamma:[0,1]\rightarrow M_\R$ is said to be generic if it does not end at walls, does not pass through the joints, and is transverse to the walls. We label the intersection points of $\gamma$ with the walls by $\gamma(t_i)$ for $t_1<t_2<\ldots$, with $\f{d}_i$ denoting the wall containing $\gamma(t_i)$. The wall-crossing operator along $\gamma$ is defined\footnote{We use the inverse limit with respect to the order $l$ for treating infinite products.} as the product:
\begin{align*}
\frp_{\gamma}=\cdots \frp_{\frd_2}^{\epsilon_2} \frp_{\frd_1}^{\epsilon_1},
\end{align*}
where $\epsilon_i\coloneqq -\sign \langle n_0(\frd_i), \gamma'(t_i)\rangle$ ($\sign$ means the sign). Two scattering diagrams $\frD_1$, $\frD_2$ are said to be equivalent if, for any generic path $\gamma$, they produce the same wall-crossing operator: $\frp^{\frD_1}_{\gamma}=\frp^{\frD_2}_{\gamma}$.

The scattering diagram $\frD$ is said to be consistent if the wall-crossing operator $\frp_{\gamma}$ is the identity for every generic loop $\gamma$. Let $\frD_{\incoming}$ denote the collection of its incoming walls. Assume all the incoming walls are hyperplanes.

\begin{thm}[\cite{gross2018canonical}]
Up to equivalence, there is a unique consistent scattering diagram whose set of incoming walls is $\frD_{\incoming}$.
\end{thm}

For any scattering diagram $\f{D}$, the positive cone $C^+=\sum_{k\in I_\ufv} \R_{\geq 0} f_k +\sum_{j\in I_\fv} \R f_j$ and the negative cone $C^-=-C^+$ will be contained in chambers of $\f{D}$. Let $\theta_{-,+}$ denote the wall-crossing operator associated to a generic path from $C^+$ to $C^-$.
\begin{thm}[\cite{gross2018canonical}]\label{thm:GHKK_determine_diagram}
Two consistent scattering diagrams $\frD_1$ and $\frD_2$ are equivalent if and only if $\theta^{\frD_1}_{-,+}=\theta^{\frD_2}_{-,+}$.
\end{thm}

\subsection*{Cluster Scattering Diagrams}

Recall that the dilogarithm function $\Li_2$ is given by
\begin{align*}
\Li_2(z)=\sum_{l\geq 1} \frac{1}{l^2}z^l.
\end{align*}
For each $k\in I_{\ufv}$, we define a wall $\frd_k=(e_k^\perp,\frp_k)$ with $\frp_k=\exp(-d_k \Li_2(-y_k))$. The cluster scattering diagram $\frD(\sd)$ associated to the initial seed $\sd$ is the consistent scattering diagram whose incoming walls are $\frd_k$, $k\in I_{\ufv}$.

The seed $\sd$ is said to be injective-reachable if there is a generic path from $C^+$ to $C^-$ which only crosses finitely many walls of $\frD(\sd)$.

Let $\sd^{\prin}$ denote the principal coefficient seed as before. As in \cite[\S 5.8]{MandelQin2021}, let $\rho$ denote the natural projection from $M(\sd^{\prin})_{\R}$ to $M_{\R}$. It sends any wall $(\frd,\frp_{\frd})$ of $\frD(\sd^{\prin})$ to a wall $(\rho(\frd),\frp_{\frd})$ in $M_{\R}$. In this way, $\rho$ gives a bijection between the walls of $\frD(\sd^{\prin})$ with those of $\frD(\sd)$.

\subsection*{Stability Scattering Diagrams}

Consider a seed $\sd$. Assume $(B_{ij})_{i,j\in I_{\ufv}}$ is skew-symmetric and let $Q$ be the corresponding quiver (without frozen vertices). Let $\C Q$ denote the path algebra of $Q$---recall that this is the algebra of $\C$-linear combinations of paths in $\Q$ (including an idempotent ``lazy path'' $e_i$ at each vertex $i$) with multiplication of paths $\alpha\beta$ given by concatenation, i.e., attaching the initial vertex of $\alpha$ to the terminal vertex of $\beta$. For any $l\in \NN$, let $\C Q^l$ denote the quotient algebra of $\C Q$ by the ideal spanned by paths of length strictly larger than $l$. The completed path algebra $\widehat{\C Q}$ is defined as the inverse limit of $\C Q^l$.  We can denote $\C Q^{\infty}= \C Q$. 

A potential $W$ of the quiver $Q$ is a $\C$-linear combination of (possibly infinitely many) closed paths in $\widehat{\C Q}$. Given an arrow $a\in Q_1$ and a closed path $w=a_1\cdots a_k$ in $\C Q$, one defines \cite{DerksenWeymanZelevinsky08}
\begin{align}\label{eq:partiala}
    \partial_a w = \sum_{i:a_i=a} a_{i+1}\cdots a_ka_1\cdots a_{i-1}.
\end{align} One extends $\partial_a$ linearly to arbitrary (possibly infinite) linear combinations of cycles in $\widehat{\C Q}$.

Assume that $B(\sd)$ is skew-symmetric and let $Q$ denote the corresponding quiver and $Q^{\op}$ its opposite quiver.  Let $W$ be a potential for $Q^{\op}$. The associated completed Jacobian algebra is
\begin{align*}
J^{\op}\coloneqq \widehat{\C Q^{\op}}/\overline{\langle \partial W\rangle },
\end{align*}
where $\partial W$ is the collection of partial derivative $\partial_a W$ for $a\in Q^{\op}_1$, $\langle \partial W\rangle$ is the corresponding two-sided ideal, and $\overline{\langle \partial W\rangle }$ is the closure of this ideal. See \S \ref{sec:covering_representation} or \cite{DerksenWeymanZelevinsky08} for details.

The stability scattering diagram $\frD^{\st}(J^{\op})$ associated to $J^{\op}$ and $\sd$ is a consistent scattering diagram encoding the following data: $\bigoplus_{k\in I_{\uf}} \R f_k$ is the space of stability conditions, and the walls are computed from the moduli spaces of finite-dimensional nilpotent semistable left $J^{\op}$-modules in the sense of \cite{king1994moduli}. We use the the motivic Hall algebras associated with nilpotent modules as in  \cite[Section 7.1.2]{Nagao10}, and refer the reader to \cite{bridgeland2017scattering} for the precise construction of the scattering diagram. For the purpose of this paper, it suffices to know the following result. 

\begin{lem}\label{lem:stability_projective}
The stability scattering diagram $\frD^{\st}(J^{\op})$ is a consistent scattering diagram such that the action of its wall-crossing operator $\theta^{\st}_{-,+}$ on $\kk\llb \cM\rrb$ satisfies
\begin{align}\label{eq:theta}
\theta^{\st}_{-,+}(x_i)=\begin{cases}
x_i\cdot (\sum_{n\in N^\oplus}\chi(\Quot_n^{\mathrm{nilp}}(P_i))x^{p^*(n)})&i\in I_{\ufv}\\
x_i&i\in I_{\fv}
\end{cases},
\end{align}
where $\Quot_n^{\mathrm{nilp}}(P_i)$ denotes the variety consisting of the $n$-dimensional nilpotent quotient modules of $P_i=J^{\op}e_i$ and $\chi$ denotes the Euler characteristic.
\end{lem}
\begin{proof} 
Let us denote the coefficients of $x_i\cdot x^{p^*(n)}$ in $\theta^{\st}_{-,+}(x_i)$ by $c_n$, and the corresponding coefficient for $\theta^{\text{Hall}}_{-,+}(x_i)$ of the motivic Hall algebra scattering diagram by $c_n^{\text{Hall}}$, see \cite[\S 6]{bridgeland2017scattering}. If $J^{\op}$ is finite-dimensional, the desired formula for $c_n$ is implied by \cite[Theorem 1.4]{bridgeland2017scattering}, see \cite[\S A.2]{qin2019bases}.

 Next, assume $J^{\op}$ is infinite-dimensional. By \cite[Theorem 10.1]{bridgeland2017scattering}, $c^{\text{Hall}}_n$ is given by the motive $[F(n,f_i,\theta)\xrightarrow{r} \mathcal{M}]$, where $\theta$ is a general point in the negative cone $C^-$, $F(n,f_i,\theta)$ a fine moduli scheme of $(n,1)$-dimensional semistable framed $J^{\op}$-representations, $\mathcal{M}$ the moduli of finite-dimensional nilpotent $J^{\op}$-representations, and the $r$ the obvious restriction map, see \cite[\S 8.1]{bridgeland2017scattering}.

By choosing $l\in \NN$ large enough, we can view these (framed) representations as (framed) representations of the truncated Jacobian algebra $(J^{\op})^l$, which is a finite-dimensional algebra. Then 
  \cite[Proposition 8.4]{bridgeland2017scattering} implies that the points of $F(n,f_i,\theta)$ are in natural bijection with those of $\Quot_n((J^{\op})^l e_i)$. By Lemma \ref{lem:quot-finite} below, $\Quot_n((J^{\op})^l e_i)$ is isomorphic to $\Quot_n^{\mathrm{nilp}}(P_i)$. Applying the integration map on $c^{\text{Hall}}_n$ (see \cite[Theorem 7.4]{Nagao10} \cite[Theorem 6.12]{joyce2007configurations}), we obtain $c_n=\chi(\Quot_n^{\mathrm{nilp}}(P_i))$ as desired.
  \end{proof}

\begin{rem}\label{rem:failure_injectivity}
As considered in \cite[Section A.2]{qin2019bases}, when $\sd$ satisfies the injectivity assumption, the action of the wall-crossing operator $\theta^{\st}_{-,+}$ is faithful. In this case, $\theta^{\st}_{-,+}$ is uniquely determined by the values in \eqref{eq:theta}. 

In general, we use $\frD^{\st,\prin}(J^{\op})$ to denote the stability scattering diagram associated to $J^{\op}$ and $\sd^{\prin}$. Then the value in \eqref{eq:theta} for $\frD^{\st}(J^{\op})$ is obtained from that for $\frD^{\st,\prin}(J^{\op})$ by evaluating $x_{i'}=1$ for all new frozen vertices $i'$. In addition, the previous projection $\rho$ also gives a bijection between the walls of $\frD^{\st,\prin}(J^{\op})$ with those of $\frD^{\st}(J^{\op})$.
\end{rem}

\begin{rem}\label{rem:dim}
    The dimension vector $n=\dim(U)\in N^{\oplus}$ of a finite-dimensional $J$- or $J^{\op}$-module $U$ is defined as follows.  For each vertex $i\in Q_0$, the component $n_i$ of $n$ is the $\C$-dimension of $e_iU$.
\end{rem}

\section{Restriction of Scattering Diagrams}\label{sec:restriction}

\subsection*{Covering of Seeds}

Let there be given equivalence relations on $I_{\ufv}$ and $I_{\fv}$ respectively. Denote the equivalence class of $i\in I$ by $\bi$ and the set of equivalent classes in $I$ by $\bI$. Define the map $\pi:I\rightarrow \bI$ such that $\pi (i)=\bi$. Let $|\bi|$ denote the cardinality of $\bi$. We make the following assumptions:
\begin{itemize}
\item $d_i=d_{i'}$ for any $i'\in \bi$.
\item For any $\bi$ and any $j'\in \bj$, we have $\sum_{i'\in \bi} B_{i',j}=\sum_{i'\in \bi} B_{i',j'}$.
\end{itemize}
Then we introduce the following data:
\begin{itemize}
\item The partition $\?{I}=\?{I}_{\ufv}\sqcup \?{I}_{\fv}\coloneqq \{\bi| i\in I_{\ufv}\}\sqcup \{\bi| i\in I_{\fv}\}$
\item The $\bI\times \bI$-matrix $\bB$ defined by $\bB_{\bi,\bj}\coloneqq \sum_{i'\in \bi} B_{i',j}$.
\item $\bd_{\bi}\coloneqq |\bi|\cdot d_i$.
\item The lattices $\bN$, $\bN_{\uf}$, $\bM$, $\bcM$, the associated basis $\?{E}=(e_{\bi})_{\bi\in\bI}$, the linear map $p^*:\bN_{\ufv}\rightarrow \bcM$, and the associated monoids such as $\bN^{\oplus}$, $\bM^{\oplus}$ are defined similarly as before. 
\end{itemize}

\begin{lem}\label{lem:project_seed}
We have $\bB_{\bi,\bj}\bd_{\bj}=-\bB_{\bj,\bi}\bd_{\bi}$.
\end{lem}
\begin{proof}
For any $j'\in \bj$, we have $\bB_{\bi,\bj}=\sum_{i'\in \bi} B_{i',j'}$. So we can rewrite
\begin{align*}
\bB_{\bi,\bj}=\frac{1}{|\bj|}\sum_{j'\in \bj}\sum_{i'\in \bi}B_{i',j'}.
\end{align*}
Similarly, we have
\begin{align*}
\bB_{\bj,\bi}=\frac{1}{|\bi|}\sum_{i'\in \bi}\sum_{j'\in \bj}B_{j',i'}.
\end{align*}
The desired claim follows from $B_{i',j'}d_{j}=-B_{j',i'}d_{i}$.
\end{proof}
By Lemma \ref{lem:project_seed}, $\?{\sd}=(\?{I}_\ufv,\?{I}_{\fv},(\bd_{\bi}),\?{E},\?{M},\bB)$ is a seed. We say $\sd$ is a covering of $\?{\sd}$ and write $\pi:\sd\rightarrow\?{\sd}$. When all orbits $\bi$ have the same cardinality $d$, we say $\pi$ is a $d:1$-covering.

Notice that, when $\sd$ is skew-symemtric, $\?{\sd}$ is skew-symmetric if and only if $\pi$ is a $d:1$ covering for some $d$.

\subsection*{Restriction}

We introduce the linear maps $\kappa:\bM_{\R}\rightarrow M_{\R}$ and $\pi:\cM\rightarrow \bcM$ such that\footnote{$\pi$ is denoted by $\iota^*$ in \cite{MandelQin2021}.}
\begin{align*}
\kappa(f_{\bi})&=\frac{1}{|\bi|}\sum_{i'\in \bi} f_{i'}\\
\pi(f_i)&=f_{\bi}
\end{align*}
Notice that we have $\pi(p^* (e_k))= p^*(e_{\bk})$ for any $k\in I_{\ufv}$.

Correspondingly, we introduce the linear map $\pi:N_{\ufv}\rightarrow \?{N}_{\ufv}$ such that $\pi(e_k)=e_{\bk}$, which induces the linear map $\pi:\hfrg\rightarrow \?{\hfrg}$ and the set-theoretic map $\pi:\hG\rightarrow \?{\hG}$.

We view $\bM_{\R}$ as a subspace of $M_{\R}$ via the embedding $\kappa$. Let $\f{D}$ denote a consistent scattering diagram in $M_\R$. Following \cite[\S A.3]{MandelQin2021}, we can restrict $\frD$ to a scattering diagram $\?{\frD}$ in $\?{M}_\R$.  More precisely, we construct $\?{\frD}$ as the inverse limit of $\?{\frD}_l$, where $\?{\frD}_l$ is defined for each $l\in \NN$ as follows: take any two generic points in $\?{M}_{\R}$ and a generic smooth path $\gamma$ connecting them (generic with respect to $\kappa^{-1} (\supp \frD_l$)). We slightly deform $\gamma$ to $\wt{\gamma}$ in $M_\R$ so that it becomes generic with respect to $\supp \frD_l$. Then $\?{\frD}_l$ is defined so that the corresponding wall-crossing operator $\?{\frp}_{\gamma}$ is equal to $\pi (\frp_{\wt{\gamma}})$.

As in \cite[A.3.3.]{MandelQin2021}, we will further assume that the equivalence classes $\bi$ are the orbits in $I$ under the action of a finite group $\Pi$ such that $\sd$ has $\Pi$-symmetry: for any $g\in \Pi$ and $i,j\in I$, we have
\begin{align*}
    B_{ij}=B_{g.i,g.j}
\end{align*}
where $g.$ denotes $g$'s action. Then $\Pi$ naturally acts on the $N^+$-graded Lie algebra $\frg$. Let $\frg^\Pi$ denote the sub Lie algebra of $\frg$ consisting of the $\Pi$-invariant elements.  The group $\Pi$ also naturally acts on $M_{\R}$, so any $g\in \Pi$ sends walls to walls.  The $\Pi$-symmetry of $\sd$ ensures that $\f{D}(\sd)$ is invariant under the $\Pi$-action.

\begin{thm}[{\cite[Theorem A.7]{MandelQin2021}}]\label{thm:consistency_restriction}
Let $\f{D}$ be a consistent scattering diagram which is invariant under the $\Pi$-action.  Then its restriction $\?{\f{D}}$ is also consistent.\footnote{\cite{MandelQin2021} stated the result for the cluster scattering diagram associated to a seed $\sd$ with $\Pi$-symmetry, but the arguments there work more generally for any $\Pi$-invariant consistent scattering diagram.}
\end{thm}

\section{Covering of Representations}\label{sec:covering_representation}
\subsection*{Covering of quivers}

Consider quivers $Q$ and $\bQ$ possibly with loops and $2$-cycles. We use $s(p)$ and $t(p)$ to denote the start point (source) and end point (target) of a path. Let $I$ and $\bI$ denote the sets of vertices for $Q$ and $\bQ$. Let $\pi$ denote a map from $Q$ to $\bQ$ respecting the quiver structures, namely:
\begin{itemize}
    \item For any $i\in I$, we have $\pi(i)\in \bI$.
    \item For any arrow $a$ from $i$ to $j$ in $Q$, $\pi(a)$ is an arrow from $\pi(i)$ to $\pi(j)$.
    \item We have $\pi(e_i)=e_{\pi(i)}$ for the lazy paths.
\end{itemize}
As in \S \ref{sec:restriction}, we denote $\pi(i)=\bi$. 
Note that the map $\pi$ induces a linear map between the vector spaces $\pi:\widehat{\C Q}\rightarrow \widehat{\C \bQ}$. 

Let $\sigma$ denote the linear map from $\widehat{\C \bQ}$ to $\widehat{\C Q}$ defined by
\begin{align*}
	\sigma(\bp)=\sum_{p\in \pi^{-1}(\bp)}p
\end{align*}	
for any path $\bp$, where $\pi^{-1}(\bp)$ denotes the set of paths in $Q$ whose image is $\bp$.
\begin{lem}\label{lem:path_alg_hom}
   The map $\sigma$ is an algebra homomorphism.
\end{lem}
\begin{proof}
    It is straightforward to reduce to the following claim: given a path $\bp$ and an arrow $\ba$ in $\bQ$, we have $$\sigma(\bp \ba)=\sigma(\bp)\sigma(\ba).$$  Every term on the left-hand side above can clearly be factored in a unique way as $pa$ for some lifts $p$ and $a$ of $\bp$ and $\ba$, respectively, giving a term on the right-hand side.  Conversely, terms on the right-hand side have the form $pa$ for lifts $p$ and $a$, and such a product is nonzero if and only if it corresponds to a term on the left-hand side. 
\end{proof}

From now on, we further assume that the covering $\pi$ is associated with a finite group $\Pi$ acting freely on $Q$ such that the following hold:
\begin{itemize}
	\item Each group element $g\in \Pi$ is a map from $Q$ to $Q$ respecting the quiver structures.
        \item The fibers $\pi^{-1}(\bk)$, $\bk\in \bQ_0$, are 
        exactly the $\Pi$-orbits in $Q_0$, and the fibers $\pi^{-1}(\ba)$, $\ba\in \bQ_1$, are exactly the $\Pi$-orbits in $Q_1$.
\end{itemize}
We denote the action of $g$ on any element $x$ by $g.x$.  Let $d=|\Pi|$.  Since $\Pi$ acts freely, each fiber has cardinality $d$.

In analogy to covering maps between topological spaces, we say $\pi:Q\rightarrow \bQ$ is a $d:1$ or $d$-folded covering map with the deck transformation group $\Pi$. When neither $Q$ nor $\bQ$ has any loops or oriented $2$-cycles, $\pi$ gives rise to a $\Pi$-symmetric $d:1$ covering between the corresponding seeds (without adding more frozen vertices) as in \S \ref{sec:restriction}.

\begin{figure}[h]
    \centering
    \begin{subfigure}[b]{0.4\textwidth}
        \centering
        \begin{tikzpicture}
            [node distance=48pt,on grid,>={Stealth[round]},bend angle=45,      pre/.style={<-,shorten <=1pt,>={Stealth[round]},semithick},    post/.style={->,shorten >=1pt,>={Stealth[round]},semithick},  unfrozen/.style= {circle,inner sep=1pt,minimum size=12pt,draw=black!100},  frozen/.style={rectangle,inner sep=1pt,minimum size=12pt,draw=black!75,fill=cyan!100},   point/.style= {circle,inner sep=1pt,minimum size=5pt,draw=black!100,fill=black!100},   boundary/.style={-,draw=cyan},   internal/.style={-},    every label/.style= {black}]

             \node[unfrozen](v1) at (0,0){$1$};    
            \node[unfrozen](v2) at (0,-48pt){$2$};
             \node[unfrozen](v3) at (48pt,0){$3$};    
            \node[unfrozen](v4) at (48pt,-48pt){$4$};

            \draw[->]  (v2) -- node[left]{$a_1$} (v1);
            \draw[->] (v2) -- node[below,pos=0.3]{$b_1$} (v3);
            \draw[->]  (v4) -- node[below,pos=0.3]{$b_2$} (v1);
            \draw[->] (v4) -- node[right]{$a_2$} (v3);

            \end{tikzpicture}
            \caption*{Quiver $Q$}
    \end{subfigure}
    \begin{subfigure}[b]{0.4\textwidth}
        \centering
\begin{tikzpicture}
[node distance=48pt,on grid,>={Stealth[round]},bend angle=45,      pre/.style={<-,shorten <=1pt,>={Stealth[round]},semithick},    post/.style={->,shorten >=1pt,>={Stealth[round]},semithick},  unfrozen/.style= {circle,inner sep=1pt,minimum size=12pt,draw=black!100},  frozen/.style={rectangle,inner sep=1pt,minimum size=12pt,draw=black!75,fill=cyan!100},   point/.style= {circle,inner sep=1pt,minimum size=5pt,draw=black!100,fill=black!100},   boundary/.style={-,draw=cyan},   internal/.style={-},    every label/.style= {black}]

 \node[unfrozen](v1) at (0,0){$\overline{1}$};    
\node[unfrozen](v2) at (0,-48pt){$\overline{2}$};

\draw[->]  (v2.120) -- node[left]{$\overline{a}$} (v1.240);
\draw[->] (v2.60) -- node[right]{$\overline{b}$} (v1.300);

\end{tikzpicture}
\caption*{Quiver $\bQ$}
    \end{subfigure}
\caption{A $2:1$ covering of a Kronecker quiver}
\label{fig:cover-quiver-kronecker}
\end{figure}
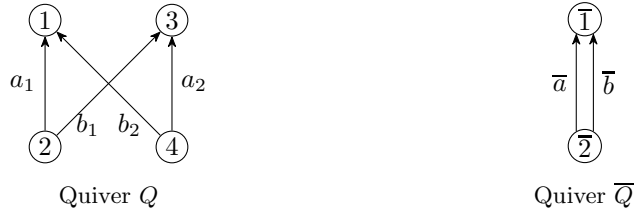

\begin{example}\label{eg:folding_quiver}
In Figure \ref{fig:cover-quiver-kronecker}, we have a $2:1$-covering map $\pi$ from $Q$ to $\bQ$ with the deck transformation group $\Pi=\{e,g\}$ such that $e$ acts trivially, and $g$ acts nontrivially on each of the sets $\{1,3\}$, $\{2,4\}$, $\{a_1,a_2\}$, and $\{b_1,b_2\}$. The map $\pi$ is given by $\pi(1)=\pi(3)=\overline{1}$, $\pi(2)=\pi(4)=\overline{2}$, $\pi(a_k)=\overline{a}$, $\pi(b_k)=\overline{b}$.
\end{example}

\begin{lem}
For any arrow $\ba$ starting from (or ending at) $\bk$, there exists exactly one arrow $a_k\in \pi^{-1}(\ba)$ starting from (resp. ending at) $k$ for each $k\in \pi^{-1}(\bk)$.
\end{lem}
\begin{proof}
Let $a$ denote any arrow in $\pi^{-1}(\ba)$. Then $a$ starts from (or ends at) some $k'\in \pi^{-1}(\bk)$. Applying the deck transformation group $\Pi$ to $a$, we obtain at least one arrow in $\pi^{-1}(\ba)$ starting from (or ending at) $k$ for each $k\in \pi^{-1}(\bk)$. The statement now follows from the fact that $\pi^{-1}(\bk)$ and the orbit $\Pi a$ have the same cardinality $d$.
\end{proof}

\begin{cor}\label{cor:lift_path}
	For any (finite or infinite) path $\bp$ starting from (or ending at) $\bk$ in $\bQ$ and any point $k\in \pi^{-1}(\bk)$, there exists exactly one path $p$ in $Q$ starting from (resp. ending at) $k$ such that $\pi(p)=\bp$.
\end{cor}
A path $p$ in $\pi^{-1}(\bp)$ is called a lift of $\bp$. Corollary \ref{cor:lift_path} has the following consequence. 
\begin{cor}\label{cor:determining_lift}
	A lift $p$ of a finite path $\bp$ is uniquely determined by its starting point, and it is also uniquely determined by its ending point. 
\end{cor}

Corollary \ref{cor:determining_lift} implies the following result.
\begin{lem}\label{lem:deck_transformation}
Let $\bi$ denote any chosen vertex of $\bQ$. If $\bQ$ is connected, then the action of a deck transformation $g$ is uniquely determined by its action on $\pi^{-1}(\bi)$. 
\end{lem}
\begin{proof}
Assume that we have determined the action of $g$ on $\pi^{-1}(\bi)$. Let $\bj$ be a vertex which is connected to $\bi$ by a path $\bp$ from $\bi$ to $\bj$. For any $i\in \pi^{-1}(\bi)$, the path uniquely lifts to a path $p$ from $i$ to some $j\in \pi^{-1}(\bj)$, and similarly, it uniquely lifts to a path $p'$ from $g.i$ to some $j'\in \pi^{-1}(\bj)$. Then the action must satisfy  $g.p=p'$, hence $g.j = j'$.  A similar argument applies to paths $\bp$ from $\bj$ to $\bi$. Repeating this progress, we obtain the desired result. 
\end{proof}

We now see that the $\Pi$-orbit of a path $p$ consists precisely of the lifts of the path $\pi(p)$.  That is,
 \begin{align}\label{eq:s_deck}
    \sigma(\pi(p))= \sum_{g\in \Pi} g.p.
\end{align}
  Take any $0\neq l\in \NN\cup\{\infty\}$.  The following result on pullbacks of certain modules follows from Corollary \ref{cor:determining_lift}.
\begin{lem}\label{lem:sigma-Pk}
    For any $k\in \pi^{-1}(\bk)$, we have $\sigma^* (\C Q^l e_k)= \C \bQ^l e_{\bk} $.
\end{lem}

\begin{rem}
We will usually view quiver representations as modules over path algebras, but let us recall that a representation $V$ of $Q$ can also be understood as a collection of vector spaces $V_k$ and linear maps $V_a: V_{s(a)}\rightarrow V_{t(a)}$ for vertices $k$ and arrows $a:s(a)\rightarrow t(a)$. The pullback by the map $\sigma$ gives a representation $\bV\coloneqq \sigma^* V$ of $\bQ$. It turns out that we have, for any $ \bk\in \bQ_0$, $\ba\in \bQ_1$,
\begin{align}
    \begin{aligned}
        \bV_{\bk}&=\bigoplus_{k\in \pi^{-1}(k)} V_k \qquad \text{and} \qquad
        \bV_{\ba}&=\bigoplus_{a\in \pi^{-1}(\ba)} V_a.
    \end{aligned}\label{eq:pullback-module}
\end{align} 
\end{rem}

\subsection*{Representations of quivers with relations}

Let there be given a collection $R$ of elements $r_i$ in $\wh{\C Q}$, indexed by $i$. Take any order $l>0$. Let $\mm$ denote the (two-sided) ideal of $\C Q^l$ generated by all arrows and assume that $r_i\in \mm^2$ for each $i$. For any $l>0$, let $\langle R^l \rangle $ denote the ideal of $\C Q^l$ generated by all $r_i$. Without loss of generality, we can assume that each $r_i$ is a finite linear combination of paths sharing the same start point and the same end point.\footnote{For any $r_i$ and vertices $j,k$, the element $e_j r_i e_k$ is contained in the ideal generated by $r_i$. So, we can replace $r_i$ by the elements $e_j r_i e_k$ for all $j,k\in Q_0$.\label{foot:ere}} We define $A^l$ to be the corresponding quotient algebra $\C Q^l/\langle R^l\rangle$. Let $\widehat{A}$ denote the inverse limit of $A^l$.

Similarly, for $\pi:Q\rar \bQ$ a $d:1$ covering, let there be given a collection $\bR$ of elements in $\C \bQ$. We construct quotient algebras $\bA^l$ and the inverse limit $\widehat{\bA}$.

Assume that $\sigma( \langle \bR^l \rangle)\subset \langle R^l\rangle $. Then we obtain an algebra homomorphism $\bsig:\bA^l\rightarrow A^l$. Correspondingly, for any $A^l$-module $V$, we have an $\bA^l$-module $\bV\coloneqq \bsig^* V$. In particular, it follows from Lemma \ref{lem:sigma-Pk} that $\bsig^*(A^l e_k)$ is a sub $\bA^l$-module of $\bA^l e_{\bk}$. Denote $P^l_k=A^l e_k$ and $\bP^l_{\bk}=\bA^l e_{\bk}$. They are indecomposable projective modules of $A^l$ and $\bA^l$ respectively.

Denote $P_k= \widehat{A} e_k$. It is the $k$-th indecomposable projective module of $\widehat{A}$. We can view it as the inverse limit of $P^l_k=A^l e_k$. Similarly define $\bP_{\bk}=\widehat{\bA}e_{\bk}=\varprojlim \bP_{\bk}^l$. Then we deduce the following result.
\begin{lem}\label{lem:folding_projective}
	$\bsig^* (P_k)$ is a sub $\widehat{\bA}$-module of $\bP_{\bk}$.
\end{lem}

Given $l\in \bb{Z}_{>0}$ and a path $p\in \C Q$, let $u_p$ denote the equivalence class of $p$ in the quotient $A^l=\C Q^l/\langle R^l\rangle$.  Similarly denote $\?{u}_{\?{p}}\in \bA^l =\C \bQ^l/\langle \bR^l\rangle$ for any path $\?{p}\in \?{Q}$.  Notice that $P^l_k$ has a basis (over $\C$) consisting of any maximal linearly independent set of elements $u_p$ for paths $p$ of length at most $l$ starting from $k$ in $Q$.  Given such a basis $\{u_p\}_p$, we have a corresponding basis $\{\pi(u_p)\}_p=\{\?{u}_{\pi(p)}\}_p$ for $\bsig^*(P^l_k)\subset \bP^l_{\bk}$ (in particular, $P_k^l$ and $\bsig^*(P_k^l)$ are isomorphic as $\C$-vector spaces).

\begin{lem}\label{lem:criterion_equality}
If $\sigma:\widehat{\C \bQ}\rightarrow \widehat{\C Q}$ descends to an injective algebra homomorphism $\bsig$ from $\bA$ to $A$, then $\bsig^* P_k =\bP_{\bk}$. 
\end{lem}
\begin{proof}
    Fix a basis $\{\?{u}_{\bp}\}_{\bp}$ for $\bP^l_{\bk}$.  Since $\bsig:\bA^l\rar A^l$ is injective and linear with image in the $\Pi$-invariant part of $A^l$, we see that $\{\bsig(\?{u}_{\bp})e_k$\} is a linearly independent subset of $P_k^l$.  Furthermore, $\pi(\bsig(\?{u}_{\bp})e_k)= \?{u}_{\bp}$, so $\{\pi(\bsig(\?{u}_{\bp})e_k)\}$ forms a basis for both $\bsig^*(P^l_k)$ and $\bP^l_{\bk}$.  Thus,  $\bsig^*(P^l_k)=\bP^l_{\bk}$ for each $l$, so $\bsig^*P_k=\bP_{\bk}$. 
\end{proof}

\subsection*{Quiver potentials and Jacobian ideals}

Let $\pi:Q\rar \bQ$ be a $d:1$ covering with deck transformation group $\Pi$.  Let $\bW$ be a potential for $\bQ$ such that $W\coloneqq \sigma(\bW)$ is a potential for $Q$.

Recall the operators $\partial_a$ as in \eqref{eq:partiala}.  Given a linear combination of arrows $\alpha=\sum_{a\in Q_1} c_a a\in \C\langle Q_1\rangle$, let $\partial_{\alpha}\coloneqq \sum_{a\in Q_1} c_a \partial_a$.

\begin{lem}\label{lem:s-partial}
 For any $\ba\in \bQ_1$ and any closed path $\?{w}$ in $\bQ$, we have
 \begin{align*}
    \sigma(\partial_{\ba}(\?{w}))=\partial_{\sigma(\ba)}\sigma(\?{w}).
 \end{align*}
\end{lem}

Let $R$ denote the collection of elements $\partial_a W$ in $\widehat{\C Q}$.  The associated completed quotient algebra $\wh{A}$ as above is called the completed Jacobian algebra and will be denoted by $J$. Similarly define $\overline{R}$ and $\bJ$. Lemma \ref{lem:s-partial} ensures that $\sigma(\partial_{\?{a}}(\bW))$ is contained in the ideal generated by $R$ for each $\?{a}\in \bQ_1$. So the map $\sigma:\widehat{\C \bQ}\rar \widehat{\C Q}$ descends to an algebra homomorphism $\bsig:\bJ\rar J$.

\begin{prop}\label{prop:s-inj}
    The map $\bsig:\bJ\rar J$ is injective.
\end{prop}
\begin{proof}
It suffices to work to arbitrary order $l>0$.
We claim that, for any $p\in \sigma(\C\bQ^l)\cap \langle R^l\rangle$, we have $p\in \sigma(\langle \bR^l\rangle)$. If so, take any $q\in \widehat{\C \bQ}$ such that $\sigma(q)$ equals $0$ in $J$. Up to any order $l$, we have that $\sigma(q)$ vanishes in $\C Q^l/\langle R^l\rangle$. Then our claim implies that $\sigma(q) \in \sigma (\langle \bR^l\rangle)$. The injectivity of $\sigma$ then implies that $q\in  \langle \bR^l\rangle$. Therefore, $q$ equals $0$ in $\bJ$.  

 Being in $\sigma(\C\bQ^l)$ ensures that $p$ is $\Pi$-invariant.  Combining this with the fact that $p\in \langle R^l\rangle$, we see that $p$ can be expressed in the form
 \begin{align}\label{eq:psum}
        p=\sum_{g\in \Pi}\sum_{i\in I_p} c_i \left(g.p_{i,1}\right)(g.\partial_{a_i}W)\left(g.p_{i,2}\right)
    \end{align}
 for some finite index-set $I_p$, coefficients $c_i\in \C$, paths $p_{i,1},p_{i,2}\in \C Q^l$, and arrows $a_i\in Q_1$.
Note that 
    \begin{align}\label{eq:sum-g-partial}
        \sum_{g\in \Pi} g.\partial_{a_i}(W)&=\partial_{\sigma(\pi(a_i))}(W)  \nonumber \\
        &= \sigma(\partial_{\pi(a_i)}(\bW)) \qquad \text{(by Lemma \ref{lem:s-partial})}.
    \end{align}
Similarly, note that 
\begin{align}\label{eq:sum-gij}
    \sum_{g\in \Pi} g.p_{i,j}= \sigma(\pi(p_{i,j}))
\end{align}
for each $i\in I_p$, $j\in \{1,2\}$.

We may assume that no term of \eqref{eq:psum} is equal to $0$ in $\C Q^l$. For any $i\in I_p$ and $j\in \{1,2\}$, since the product $\left(g.p_{i,1}\right)(g.\partial_{a_i}W)\left(g.p_{i,2}\right)$ is nonzero for each $g\in \Pi$, we see using Corollary \ref{cor:determining_lift} that, for any $g,g',g''\in \Pi$, the product $\left(g''.p_{i,1}\right)(g'.\partial_{a_i}W)\left(g.p_{i,2}\right)$ equals $0$ unless $g'.\partial_{a_i}W=g.\partial_{a_i}W$ and $g''.p_{i,1}=g.p_{i,1}$. 
We may thus rewrite \eqref{eq:psum} as
    \begin{align*}
        p&=\sum_{i\in I_p} c_i \left(\sum_{g''\in \Pi} g''.p_{i,1}\right)\left( \sum_{g'\in \Pi} g'.\partial_{a_i} W \right)\left(\sum_{g\in \Pi} g\cdot p_{i,2}\right)\\
        &=\sum_{i\in I_p}c_i \sigma(\pi(p_{i,1}))\sigma \left(\partial_{\pi(a_i)} \bW\right) \sigma(\pi(p_{i,2})) \qquad \text{(applying \eqref{eq:sum-g-partial} and \eqref{eq:sum-gij})}\\
        &=\sigma \left(\sum_{i\in I_p}c_i\pi(p_{i,1})\left(\partial_{\pi(a_i)} \bW\right)\pi(p_{i,2})\right) \\
        &\in \sigma(\langle \bR^l\rangle) 
        \end{align*}
    as desired.   
\end{proof}

We obtain the following as an immediate corollary of Lemma \ref{lem:criterion_equality} combined with Proposition \ref{prop:s-inj}.
\begin{prop}\label{prop:sigma-Pk}
For $P_k=Je_k$, $\bP_{\bk}=\bJ e_{\pi(k)}$, and $\bsig:\bJ\rar J$ as above, we have $\bsig^* P_k = \bP_{\bk}$.    
\end{prop}

\subsection*{Euler Characteristics of Quiver Grassmannians}
Recall that, for any $J$-module $V$ and $n\in \NN^{Q_0}$, when $V$ is finite-dimensional, the quiver Grassmannian $\Gr_n(V)$ is the projective variety parameterizing the $n$-dimensional submodules of $V$. Quiver Grassmannians for finite-dimensional $\bJ$-modules are defined similarly. Additionally, for each $l\in \NN$, we can consider $\Gr_n(V^l)$, the projective variety parametrizing $n$-dimensional submodules of $V^l=J^l\otimes V$ over $J^l$.

 Let $\pi:Q\rightarrow \bQ$ denote a covering as before.  As in \cite{haupt2012euler}, we will define $\C^*$-actions on our quiver Grassmannians by introducing $\Z$-gradings on the modules.  We will then use the following lemma to relate the stability scattering diagram of $Q$ to that of $\bQ$.

\begin{lem}[\cite{bialynicki1973some}]\label{lem:fixed_pt_Euler}
Let $\C^*$ act on a nonempty locally closed subset $X$ of a projective variety $Y$. Denote by $X^{\mathbb{C}^*}$ the fixed point subset. Then $X^{\mathbb{C}^*}$ is also a nonempty locally closed subset of $Y$, and $\chi(X)=\chi(X^{\mathbb{C}^*})$, where $\chi$ denotes the Euler characteristic.
\end{lem}

    Recall that, at least up to any finite order $l\in \bb{Z}_{>0}$, the $k$-th indecomposable projective module $J^le_k$ of $J^l$ has a basis consisting of a maximal linearly independent set of elements $u_p$ represented by paths $p$ starting at the vertex $k$.  We take the inverse limit with respect to $l$ to obtain a fixed choice\footnote{When $P_k$ is infinite-dimensional, the set of elements $u_p$ is constructed as an inverse limit with respect to the inclusion partial order.} of topological basis $\{u_p\}_p$ for $P_k=Je_k$ (an ordinary basis if $P_k$ is finite-dimensional), and a topological basis $\{\bu_{\bp}\}$ for $\bP_{\bk}$ with $\bu_{\bp}\coloneqq \pi(u_p)$. 
   
    Let $\supp_{Q_0} (P_k)$ denote the set of vertices $i$ such that $(P_k)_i\neq 0$; i.e., such that $i$ is the endpoint of a path $p$ which represents a nonzero element of $P_k$.  Let $\supp_{Q_1}(P_k)$ be the set of arrows $a$ such that $(P_k)_a\neq 0$; i.e., such that $a$ is an arrow in a path $p$ whose action on $P_k$ is non-zero. 
    
    Let $\wh{\C Q}_k$ and $J_k$ denote the quotients of $\wh{\C Q}$ and $J$, respectively, by the two-sided ideals generated by arrows $a\in Q_1\setminus \supp_{Q_1}(P_k)$.  We similarly define $\supp_{\bQ_1}(\bP_{\bk})$, $\wh{\C \bQ}_{\bk}$, and $\bJ_{\bk}$ .  Note then that we may naturally view $P_k$ as a $J_k$-module and $\bP_{\bk}$ as a $\bJ_{\bk}$-module (because the ideal we mod out by is contained in the left annihilator of $P_k$ or $\bP_{\bk}$), and that the injective homomorphism $\bsig:\bJ\rightarrow J$ descends to an injective homomorphism $\bsig:\bJ_{\bk}\rar J_k$ by Proposition \ref{prop:sigma-Pk}. We make the following assumption.

    \begin{assumption*}[Nice grading]
        There exists a $\Z$-grading $\partial$ on $\supp_{Q_0} (P_k)$ such that the following properties are true:
        \begin{itemize}
            \item For $i,i'\in \supp_{Q_0}(P_k)$ in the same fiber $\pi^{-1}(\bi)$, we have $\partial (i)=\partial(i')$ if and only if $i=i'$.
            \item For any arrows $a,a'\in \supp_{Q_1}(P_k)$ such that $\pi(a)=\pi(a')$, 
            we have $\partial(t(a))-\partial(s(a))=\partial(t(a'))-\partial(s(a'))$.
        \end{itemize}
    \end{assumption*}

    Such a grading is called a nice grading on $\supp_{Q_0}P_k$ with respect to $P_k$. We will also say $P_k$ has a nice grading for simplicity. In this case, we define a $\Z$-grading $\partial$ on $\supp_{Q_1}(P_k)$ by $\partial(a)\coloneqq \partial(t(a))-\partial(s(a))$.  We further define $\partial(\pi(a))=\partial(a)$ for $a\in \supp_{Q_1}(P_k)$, which is independent of the choice of $a$ by the second condition for $\partial$ being a nice grading.  Now $\partial$ determines gradings on $\wh{\C Q}_k$ and $\wh{\C \bQ}_{\bk}$.  Furthermore, each Jacobian ideal generator $\partial_a W$ which does not vanish in $\wh{\C Q}_k$ will be homogeneous with respect to the grading on $\wh{\C Q}_k$, so the grading descends to $J_k$.  Similarly, the grading on $\wh{\C \bQ}_{\bk}$ descends to a grading on $\bJ_{\bk}$.  By construction, $\bsig:\bJ_{\bk}\rar J_k$ respects these gradings.

    \begin{example}\label{eg:nice_grading_kronecker}
        Let us continue Example \ref{eg:folding_quiver}. Consider the projective module $P_2$ for $A=\C Q$. It has the basis $e_2,a_1,b_1$. So $\supp_{Q_0}(P_2)=\{2,1,3\}$ and $\supp_{Q_1}(P_2)=\{a_1,b_1\}$. We assign grading $\partial(2)=\partial(1)=0$ and $\partial (3)=1$. Correspondingly, $\partial(a_1)=0$ and $\partial(b_1)=1$. Then $\partial$ is a nice grading.

        Note that we can NOT extend $\partial$ to a $\Z$-grading on $Q_0$ such that, for $\partial(b_2)\coloneqq \partial(1)-\partial(4)$, $\partial(a_2)\coloneqq \partial(3)-\partial(4)$, we have $\partial(a_2)=\partial(a_1)$ and $\partial(b_2)=\partial(b_1)$.
    \end{example}
   
   We define $\C^*$-actions on $J_k$ and $\bJ_{\bk}$ (hence also on their ideals) via $t.x=t^{\partial(x)}x$ for each homogeneous element $x$. Take any $l\in\NN$. Since the induced actions on $P^l_k$ and $\bP^l_{\bk}$ take submodules to submodules and preserve dimension vectors, they induce $\C^*$-actions on the quiver Grassmannians $\Gr_n(P^l_k)$ and $\Gr_{\bn}(\bP^l_{\bk})$ for each $n\in \NN^{Q_0}$ and $\bn\in \NN^{\bQ_0}$.

\begin{rem}\label{rmk:nice-Pl}
       For each $l\in \NN$, one can define $\supp_{Q_0} (P_k^l)$ and $\supp_{Q_1} (P_k^l)$ similarly to above by using $P_k^l$ and $J^l$ in place of $P_k$ and $J$.  One can then ask for the existence of a nice grading on $\supp_{Q_0}(P_k^l)$, defined in the obvious way.  This is indeed sufficient for the proofs of Proposition \ref{prop:projection_Euler} and Corollary \ref{cor:projection_Euler} below  via essentially the same arguments (assuming $l$ in Proposition \ref{prop:projection_Euler} is this $l$, and assuming the $\bn$ in Corollary \ref{cor:projection_Euler} satisfies $|\bn|\leq l+1$).  We shall use this in \S \ref{sec:gen}.
   \end{rem}

\begin{prop}\label{prop:projection_Euler}
Take any $k\in Q_0$ and $l\in \NN$. Then we have $\chi(\Gr_{\bn}(\bP_{\bk}^l))=\sum_{n\in \pi^{-1}(\bn)}\chi(\Gr_n(P_k^l))$ for $\bn\in \NN^{\bQ_0}$.
\end{prop}
\begin{proof}
Recall from Proposition \ref{prop:sigma-Pk} that $\bsig^* P_k=P_{\bk}$.  Note that the same holds when working over $J_k$ and $\bJ_{\bk}$, as we do here.  Furthermore, since $\pi$ preserves the lengths of paths, this descends to $\bsig^* P_k^l=P_{\bk}^l$.

    Given a submodule $U$ of $P_k^l$ with $n=\dim U$, $\bsig^* U$ is a submodule of $\bP_{\bk}^l$ with $\dim \bsig^* U=\bn\coloneqq \pi(n)$.  This gives an embedding $\bsig^*:\Gr_n (P_k^l)\hookrightarrow \Gr_{\bn}(\bP_{\bk}^l)$ respecting the $\C^*$-actions.

    Let $\bU$ be an $\bn$-dimensional submodule of $\bP_{\bk}^l$.  In particular, $\bU$ has a basis consisting of linear combinations $\?{\beta}_j=\sum c_{\bp} \bu_{\bp}$ for paths $\bu_{\bp}$ starting at $\bk$.  Now suppose that $\bU$ is fixed by the $\C^*$-action.  Then $\bU$ is homogeneous; i.e., we may assume that each $\?{\beta}_j$ is homogeneous.\footnote{To see this, take a generic $t\in \C^*$. From $t.\bU=\bU$, one can deduce that $\bU$ is spanned by eigenvectors of $t$.}  Furthermore, by an argument like that of Footnote \ref{foot:ere}, we may assume for each $j$ that all paths $\bp$ contributing to $\?{\beta}_j$ end at a common vertex $\bv_j$. Then each $\bp$ is $\pi(p)$ for some path $p$ from $k$ to a vertex $v\in \pi^{-1}(\bv_j)$. By the homogeneity of $\?{\beta}_j$ and the first condition of a nice grading, $v$ must actually be the same for each such $p$; call this common vertex $v_j$.  Hence, the lift $\beta_j\coloneqq \sum c_p u_p$ of $\?{\beta}_k$ is homogeneous with respect to the dimension vector grading $\dim$ (cf. Remark \ref{rem:dim}), and $\pi(\dim \beta_j)=\dim(\?{\beta}_j)=e_{\bv_j}$.

    Now let $U$ be the $\C$-span of $\{\beta_j\}_j$.  Note that $\beta_j=\bsig(\?{\beta}_j)e_k$, so $U=\bsig(\bJ_{\bk}^l\langle \?{\beta}_j\rangle_j)e_k$, where $\bJ_{\bk}^l\langle \?{\beta}_j\rangle_j$ is an ideal generated by all $\?{\beta_j}$.  We claim that $U$ is a sub $J_k^l$-module of $P_k^l=J_k^l e_k$.\footnote{This claim can also be verified by analyzing vector space structures as in \eqref{eq:pullback-module}.}  To see this, we must check that $p\beta_j\in U$ for all $p\in J_k^l$. It suffices to consider where $p$ is represented by a single path starting at the final vertex $v$ of the paths in $\beta_j$. 
 Then Corollary \ref{cor:determining_lift} ensures that $\bsig(\pi(p))\beta_j=p\beta_j$. Combining with $\beta_j=\bsig(\?{\beta}_j)e_k$, we have $$p\beta_j=p\beta_j e_k=\bsig(\pi(p))\bsig(\?{\beta}_j)e_k=\bsig(\pi(p)\?{\beta}_{j})e_k\in \bsig(\bJ_{\bk}^l\langle \?{\beta}_j\rangle_j)e_k=U.$$  It follows from $\pi(\dim \beta_j)=\dim(\?{\beta}_j)$ that $\pi(\dim(U))=\dim(\bU)$.  Furthermore, we see that $\bsig^* U = \bU$.

So the $\C^*$-fixed point subset in $\Gr_{\bn}(\bP_{\bk}^l)$ has the following decomposition:
\begin{align*}
(\Gr_{\bn}(\bP_{\bk}^l))^{\C^*}= \bigsqcup_{n\in \pi^{-1}(\bn)} (\bsig^* \Gr_n(P_k^l))^{\C^*} \simeq \bigsqcup_{n\in \pi^{-1}(\bn)}(\Gr_n(P_k^l))^{\C^*}.
\end{align*}
The desired result follows from Lemma \ref{lem:fixed_pt_Euler}.
\end{proof}

\begin{example}
    Let us continue Example \ref{eg:nice_grading_kronecker}. Notice that $\bsig^*(P_2)=\bP_{\overline{2}}$.  Take dimension vector $\bn$ such that $(\bn)_{\overline{1}}=1$, $(\bn)_{\overline{2}}=0$. Then the quiver Grassmannian $\Gr_{\bn}(\bsig^*(P_2))$ is $\mathbb{P}^1$ and so has Euler characteristic $2$. Its $\C^*$-fixed locus consists of two points $\{\C u_{\ba},\C u_{\?{b}}\}$ and is isomorphic to the disjoint union of the quiver Grassmannians $\Gr_{(1,0,0,0)}(P_2)=\{\C u_{a_1}\}$ and $\Gr_{(0,0,1,0)}(P_2)=\{\C u_{b_1}\}$.
\end{example}

The following is essentially \cite[Lemma 6.2, Corollary 6.3]{Nagao10} but with a stronger bound.

\begin{lem}\label{lem:nilpotent-mod}
    Fix any $n\in N^{\oplus}$ and $l\in \NN$ with $l\geq |n|-1$.  Then, for any $n$-dimensional nilpotent $J$-module $V$, we have $u_p V=0$ whenever $p$ is a path of length $l+1$. In particular, the moduli stack of the $n$-dimensional nilpotent $J$-modules is canonically isomorphic to that of $n$-dimensional $J^l$-modules.
\end{lem}
\begin{proof}
    Let $\mm$ denote the two-sided ideal in $J$ generated by $u_a$ for $a\in Q_1$. Denote $V^{(d)}=\mm^d(V)$ for $d\in \NN$. We want to show that $V^{(l+1)}=0$. It suffices to verify $V^{(|n|)}=0$.

Note that we have a chain of subspaces $V=V^{(0)}\supset V^{(1)}\supset \cdots \supset V^{(d)}\supset V^{(d+1)}\supset \cdots$. If $V^{(d+1)}=V^{(d)}\neq 0$ for some $d\in\NN$, then $V^{(d')}=V^{(d)}\neq 0$ for all $d'\geq d$, contradicting to the assumption that $V$ is nilpotent. Therefore, $\dim V^{(d+1)}<\dim V^{(d)}$ whenever $V^{(d)}\neq 0$. We deduce that $V^{(|n|)}=0$.
\end{proof}

\begin{lem}\label{lem:quot-finite}
    Fix $n\in N^{\oplus}$ and $l\in \NN$ with $|n|\leq l+1$.  Then $\Quot_n(P^l_k)=\Quot_n^{\mathrm{nilp}}(P_k)$.
\end{lem}
\begin{proof}
    Let $E$ denote any $n$-dimensional nilpotent $J$-module, which is also a $J^l$-module by Lemma \ref{lem:nilpotent-mod}. Let $f$ denote any surjective homomorphism from $P_k$ to $E$. The equivalence classes $u_p$ in $P_k$ are sent to $u_p f(e_k)$, which vanish when the paths $p$ have lengths of at least $l$, see Lemma \ref{lem:nilpotent-mod}. Therefore, $f$ can be viewed as a surjective homomorphism from $P_k^l$ to $E$. Conversely, any surjective homomorphism $g$ from $P_k^l$ to $E$ induces a surjective homomorphism from $P_k$ to $E$ sending $e_k$ to $g(e_k)$. 
\end{proof}

\begin{cor}\label{cor:projection_Euler}
For each $\bn\in \NN^{\bQ_0}$, we have $\chi(\Quot_{\bn}^{\mathrm{nilp}}(\bP_{\bk}))=\sum_{n\in \pi^{-1}(\bn)}\chi(\Quot_n^{\mathrm{nilp}}(P_k))$.
\end{cor}
\begin{proof}
Fix $l\geq|\bn|-1$.  Notice that we have isomorphisms $\Quot_n^{\mathrm{nilp}}(P_k)\cong \Quot_n (P_k^l)$ (by Lemma \ref{lem:quot-finite}) and $\Quot_n (P_k^l)\cong  \Gr_{\dim P_k^l-n} (P_k^l)$. Similarly for $\bP_{\bk}$ and $\bP_{\bk}^l$.  We deduce that
\begin{align*}
    \chi(\Quot_{\bn}^{\mathrm{nilp}}(\bP_{\bk}))&=\chi(\Quot_{\bn}(\bP_{\bk}^l)) \\
    &=\chi(\Gr_{\dim \bP_{\bk}^l-\bn}(\bP_{\bk}^l))\\
                                &=\sum_{n': \pi(n')=\dim \bP_{\bk}^l-\bn}\chi(\Gr_{n'}(P_k^l)) \qquad \text{(Proposition \ref{prop:projection_Euler})}\\
                                &=\sum_{n: \pi(n)=\bn}\chi(\Gr_{\dim P_k^l -n}(P_k^l)) \qquad \text{(set $n\coloneqq\dim P_k^l -n'$)}\\
                                &=\sum_{n:\pi(n) =\bn}\chi(\Quot_n(P_k^l)) \\
                                &=\sum_{n:\pi(n)=\bn}\chi(\Quot_n^{\mathrm{nilp}}(P_k)).
\end{align*}
\end{proof}

\begin{proof}[Proof of Theorem \ref{thm:stability_diagram}]
    It suffices to verify the theorem for principal coefficients, see Remark \ref{rem:failure_injectivity}. In this case, the desired claim follows from Lemma \ref{lem:stability_projective} and Corollary \ref{cor:projection_Euler}. 

    More precisely, consider the principal coefficients stability scattering diagrams $\frD_1\coloneqq \frD^{\st,\prin}(J^{\op})$ and $\frD_2\coloneqq \frD^{\st,\prin}(\bJ^{\op})$. By Theorem \ref{thm:GHKK_determine_diagram}, they are determined by their wall-crossing operators $\theta^{\frD_i}_{-,+}$, $i=1,2$.
    Similarly, the scattering diagram $\frD_3\coloneqq \?{\frD^{\st,\prin}(J^{\op})}$ is determined by its wall-crossing operator $\theta^{\frD_3}_{-,+}$. 
     
      Note that $\theta^{\frD_i}_{-,+}$, $i=1,2$, appear as $\theta^{\st}_{-,+}$ in Lemma \ref{lem:stability_projective}. For any $k\in I$, using Corollary \ref{cor:projection_Euler}, we compute that  $\pi(\theta^{\frD_1}_{-,+}(x_{k}))=\theta^{\frD_2}_{-,+}(x_{\bk})$, where $\pi$ denotes the linear map sending $x_i$ to $x_{\bi}$ for any $i$. Moreover, we have $\theta^{\frD_3}_{-,+}(x_{\bk})=\pi(\theta^{\frD_1}_{-,+}(x_{k}))$---by construction, if $\gamma$ is a path in $\?{M}_{\bb{R}}$ from $\?{C}^+$ to $\?{C}^-$ and $\wt{\gamma}$ is a deformation of $\gamma$ giving a path from $C^+$ to $C^-$ in $M_{\bb{R}}$, then $\?{\f{p}}_{\gamma}=\pi(\f{p}_{\wt{\gamma}})$. Therefore, $\theta^{\frD_2}_{-,+}(x_{\bk})=\theta^{\frD_3}_{-,+}(x_{\bk})$ for all $\bk$. Since we work with principal coefficients, the actions of $\theta^{\frD_i}_{-,+}$, $i=1,2,3$, are faithful. Therefore, $\frD_3$ and $\frD_2$ are equivalent by Theorem \ref{thm:GHKK_determine_diagram}. The equivalence extends to arbitrary coefficients by Remark \ref{rem:failure_injectivity}.
\end{proof}

\section{Applications to Surfaces}

\subsection{Projection of Surface Jacobian Algebras}

\begin{figure}[htb]
    \centering
\tikzset{every picture/.style={line width=0.75pt}} 
\begin{tikzpicture}[x=0.75pt,y=0.75pt,yscale=-1,xscale=1]

\draw   (279,124.75) .. controls (279,93.96) and (314.37,69) .. (358,69) .. controls (401.63,69) and (437,93.96) .. (437,124.75) .. controls (437,155.54) and (401.63,180.5) .. (358,180.5) .. controls (314.37,180.5) and (279,155.54) .. (279,124.75) -- cycle ;
\draw [color={rgb, 255:red, 65; green, 117; blue, 5 }  ,draw opacity=1 ]   (351.33,179.99) .. controls (346.4,173.97) and (344.89,165.82) .. (345.4,158.1) .. controls (346.17,146.39) and (351.56,135.65) .. (356.61,134.74) ;
\draw [color={rgb, 255:red, 65; green, 117; blue, 5 }  ,draw opacity=1 ] [dash pattern={on 0.84pt off 2.51pt}]  (356.61,134.74) .. controls (372.75,138.28) and (378.65,174.94) .. (351.33,179.99) ;
\draw [color={rgb, 255:red, 65; green, 117; blue, 5 }  ,draw opacity=1 ]   (345.4,158.1) .. controls (359.63,167.79) and (384.7,170.09) .. (413.57,164.52) ;
\draw [color={rgb, 255:red, 65; green, 117; blue, 5 }  ,draw opacity=1 ]   (337.3,129.31) .. controls (326.4,136.49) and (329.4,148.78) .. (345.4,158.1) ;
\draw [color={rgb, 255:red, 65; green, 117; blue, 5 }  ,draw opacity=1 ] [dash pattern={on 0.84pt off 2.51pt}]  (377.62,125.05) .. controls (408.67,129.09) and (421.64,148.09) .. (413.57,164.52) ;
\draw   (356.61,134.74) .. controls (383.61,133.88) and (389.06,105.19) .. (353.64,106.49) .. controls (318.23,107.8) and (329.61,135.6) .. (356.61,134.74) -- cycle ;
\draw  [color={rgb, 255:red, 65; green, 117; blue, 5 }  ,draw opacity=1 ] (345.4,158.1) .. controls (408.12,158.66) and (442.99,86.93) .. (363.99,85.63) .. controls (284.99,84.32) and (282.68,157.54) .. (345.4,158.1) -- cycle ;

\draw (331.66,157.25) node [anchor=north west][inner sep=0.75pt]    {$p$};

\end{tikzpicture}
    \caption{An ideal triangulated torus with one puncture $p$}
    \label{fig:ideal-triangulated-torus}
\end{figure}
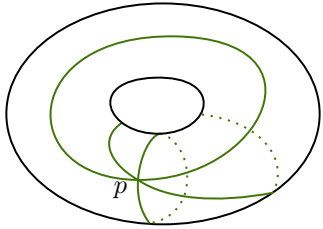

	Let $S$ be a compact oriented closed surface with punctures $M$ such that $(S,M)$ is triangulable; i.e., $\#M \geq 1$, and if $S$ is a sphere then $\#M\geq 4$. Fix an ideal triangulation $T$ of $\Sigma=(S,M)$; that is, $T$ is a maximal collection of pairwise non-intersecting (except at their endpoints) non-homotopic simple arcs between punctures.  See Figure \ref{fig:ideal-triangulated-torus}. We assume that $T$ does not cut out any self-folded triangles (i.e., triangles with a repeated edge). We refer the reader to \cite[\S 3]{MandelQin2021} for details.

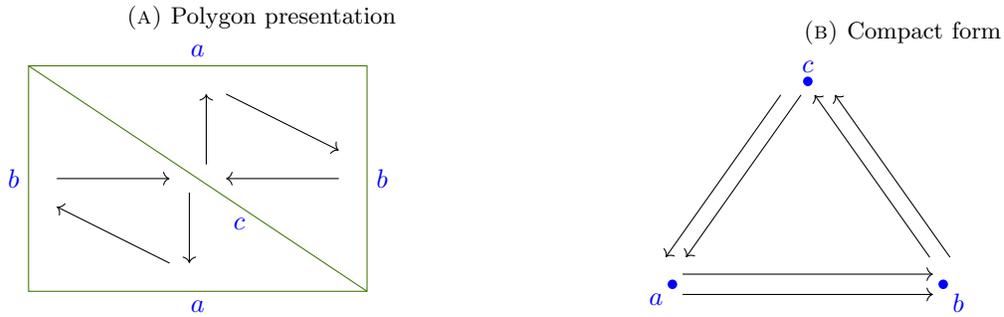
\begin{figure}[htbp]
    \centering
    \begin{subfigure}[b]{0.45\textwidth}
    \caption{Polygon presentation}
     \begin{tikzpicture}[scale=0.75]
	\draw [color={rgb, 255:red, 65; green, 117; blue, 5 }] (0,0) rectangle (6,4);
	\draw [color={rgb, 255:red, 65; green, 117; blue, 5 } ]   (0,4) -- (6,0) ;
	\fill [color=blue](3,0)node[below]{$a$};\fill(3,4)[color=blue]node[above]{$a$};\fill[color=blue](0,2)node[left]{$b$};\fill[color=blue](6,2)node[right]{$b$};
	\fill [color=blue] (4,1.2)node[left]{$c$};
	\draw [<-] (3.5,2) -- (5.5,2);	\draw [<-] (3.15,3.5) -- (3.15,2.25);	\draw [<-]  (5.5,2.5) -- (3.5,3.5);
	\draw [<-] (2.5,2) -- (0.5,2); \draw[<-](2.85,0.5) -- (2.85,1.75);\draw [<-] (0.5,1.5)--(2.5,0.5);
    \end{tikzpicture}
    \end{subfigure}
    \hfill
    \begin{subfigure}[b]{0.45\textwidth}
    \caption{Compact form} 
    \begin{tikzpicture}[scale=0.9]
	\fill(0,0)[color=blue]circle(2pt)node[below left]{$a$};	\fill[color=blue](4,0)circle(2pt)node[below right]{$b$};	\fill[color=blue](2,3)circle(2pt)node[above]{$c$};
	\draw[->](0.15,0.15)--(3.85,0.15);\draw[->](0.15,-0.15)--(3.85,-0.15); \draw[->](3.8,0.4)--(2.1,2.8); \draw[->](4.1,0.4)--(2.4,2.8); \draw[<-](0.2,0.4)--(1.9,2.8); \draw[<-](-0.1,0.4)--(1.6,2.8); 
    \end{tikzpicture}
\end{subfigure}
    \caption{The adjacency quiver for the torus, shown in polygon presentation and compact form}
    \label{fig:adjacency-quiver-torus}
\end{figure}

 Next we define the adjacency quiver $Q(T)$.  The vertices of $Q(T)$ are the arcs from triangulation $T$. The arrows in $Q(T)$ are defined as follows: There is an arrow $i\to j$ whenever arc $i$ and $j$ have a common endpoint in $M\subset S$ and $i$ is followed immediately counterclockwise by $j$ around the common endpoint.  See Figure \ref{fig:adjacency-quiver-torus}.

\begin{figure}[h]
    \centering

\tikzset{every picture/.style={line width=0.75pt}} 

\begin{tikzpicture}[x=0.75pt,y=0.75pt,yscale=-1,xscale=1]

\draw   (363.79,148.44) -- (194.17,149.02) -- (261.42,44.5) -- cycle ;
\draw    (265.9,141.13) -- (235.08,105.49) ;
\draw [shift={(233.77,103.98)}, rotate = 49.14] [color={rgb, 255:red, 0; green, 0; blue, 0 }  ][line width=0.75]    (10.93,-3.29) .. controls (6.95,-1.4) and (3.31,-0.3) .. (0,0) .. controls (3.31,0.3) and (6.95,1.4) .. (10.93,3.29)   ;
\draw    (194.17,149.02) -- (206.13,221.97) ;
\draw    (119.45,91.27) -- (194.17,149.02) ;
\draw    (194.17,149.02) -- (108.24,194.38) ;
\draw    (212.85,96.76) -- (149.84,97.35) ;
\draw [shift={(147.84,97.37)}, rotate = 359.46] [color={rgb, 255:red, 0; green, 0; blue, 0 }  ][line width=0.75]    (10.93,-3.29) .. controls (6.95,-1.4) and (3.31,-0.3) .. (0,0) .. controls (3.31,0.3) and (6.95,1.4) .. (10.93,3.29)   ;
\draw    (144.85,119.86) -- (139.86,161.21) ;
\draw [shift={(139.62,163.19)}, rotate = 276.88] [color={rgb, 255:red, 0; green, 0; blue, 0 }  ][line width=0.75]    (10.93,-3.29) .. controls (6.95,-1.4) and (3.31,-0.3) .. (0,0) .. controls (3.31,0.3) and (6.95,1.4) .. (10.93,3.29)   ;
\draw    (363.79,148.44) -- (316.71,231.33) ;
\draw    (283.09,155.11) -- (325.62,189.65) ;
\draw [shift={(327.17,190.91)}, rotate = 219.08] [color={rgb, 255:red, 0; green, 0; blue, 0 }  ][line width=0.75]    (10.93,-3.29) .. controls (6.95,-1.4) and (3.31,-0.3) .. (0,0) .. controls (3.31,0.3) and (6.95,1.4) .. (10.93,3.29)   ;
\draw    (240.5,97.24) -- (294.54,97.24) ;
\draw [shift={(296.54,97.24)}, rotate = 180] [color={rgb, 255:red, 0; green, 0; blue, 0 }  ][line width=0.75]    (10.93,-3.29) .. controls (6.95,-1.4) and (3.31,-0.3) .. (0,0) .. controls (3.31,0.3) and (6.95,1.4) .. (10.93,3.29)   ;
\draw    (304.76,103.69) -- (277.58,138.96) ;
\draw [shift={(276.36,140.55)}, rotate = 307.61] [color={rgb, 255:red, 0; green, 0; blue, 0 }  ][line width=0.75]    (10.93,-3.29) .. controls (6.95,-1.4) and (3.31,-0.3) .. (0,0) .. controls (3.31,0.3) and (6.95,1.4) .. (10.93,3.29)   ;
\draw    (210,193.33) -- (257.57,157.91) ;
\draw [shift={(259.18,156.72)}, rotate = 143.33] [color={rgb, 255:red, 0; green, 0; blue, 0 }  ][line width=0.75]    (10.93,-3.29) .. controls (6.95,-1.4) and (3.31,-0.3) .. (0,0) .. controls (3.31,0.3) and (6.95,1.4) .. (10.93,3.29)   ;
\draw    (363.79,148.44) -- (440,148.92) ;

\draw (239.73,123.23) node [anchor=north west][inner sep=0.75pt]    {$\alpha $};
\draw (312.08,162.49) node [anchor=north west][inner sep=0.75pt]    {$\overline{\alpha }$};
\draw (252.14,79.92) node [anchor=north west][inner sep=0.75pt]    {$f( \alpha )$};
\draw (287,120.66) node [anchor=north west][inner sep=0.75pt]    {$f^{2}( \alpha )$};
\draw (166.34,78.19) node [anchor=north west][inner sep=0.75pt]    {$g( \alpha )$};
\draw (101.82,132.26) node [anchor=north west][inner sep=0.75pt]    {$g^{2}( \alpha )$};
\draw (160.59,173.82) node [anchor=north west][inner sep=0.75pt]    {$\cdots $};
\draw (220.59,181.87) node [anchor=north west][inner sep=0.75pt]    {$g^{n_{\alpha } -1}( \alpha )$};

\end{tikzpicture}

    \caption{Maps on arrows}
    \label{fig:maps-on-arrows}
\end{figure}
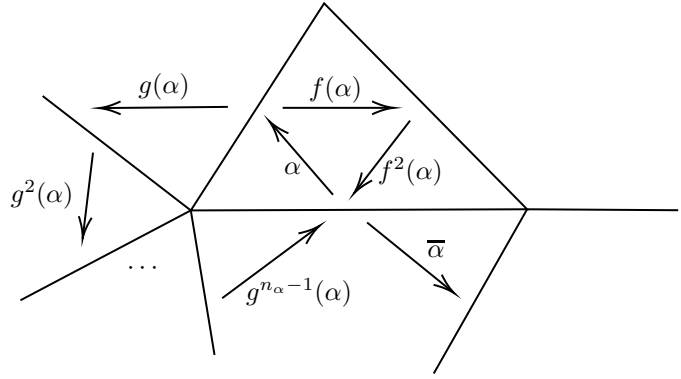

We next introduce maps defined on arrows of $Q(T)$ as introduced in \cite{Ladkani2012-zx}. Let $\alpha$ be the arrow $i\to j$ in $T$.  Then $i$ and $j$ have a common endpoint in $M$, and $i$ is followed immediately counterclockwise by $j$. There exists a unique arc $k$ sharing the same endpoint with $i,j$, which follows immediately counterclockwise of $j$. Then $g(\alpha)$ is the arrow $j\to k$. On the other hand, by definition of triangulation (without self-folded triangles), there exists some arc $l$ such that $l,i,j$ enclose a triangle of $T$. Then $f(\alpha)$ is the arrow $j\to l$. $f(\alpha)$ and $g(\alpha)$ are the two arrows in $Q(T)$ starting at $j$. It is also shown that for any vertex $i$ there exist exactly two arrows starting and two arrows ending at $i$ in $Q(T)$. Thus for an arrow $\alpha$, we can denote the other arrow sharing the same start with $\alpha$ by $\overline{\alpha}$.  See Figure \ref{fig:maps-on-arrows}.

It is easy to see that $\alpha,f(\alpha),f^2(\alpha)$ form a 3-cycle path, and $\alpha,g(\alpha),g^2(\alpha),\cdots$ rotate around some puncture. Let $n_\alpha\in\mathbb{N}$ be the least positive integer such that $g^{n_\alpha}(\alpha)=\alpha$. Then Labardini-Fragoso's association goes as follows: 

Associate to every puncture $m_i\in M$ a complex number $c_i\ne0$.  Then we define the potential (see \cite{Lad09,Ladkani2012-zx})
\begin{align}\label{eq:surface_potential}
   W=\sum f^2(\alpha) \cdot f(\alpha) \cdot \alpha-\sum c_\beta  g^{n_\beta-1}(\beta) \cdots g(\beta) \cdot \beta. 
\end{align}
Here, the first sum is over all $f$-orbits, with $\alpha$ being an arbitrary representative of the $f$-orbit.  The second sum is over all $g$-orbits, with $\beta$ an arbitrary choice of representative of the $g$-orbit. The number $c_\beta$ equals the complex number associated to the common endpoint of arcs $s(\beta),s(g(\beta)),\cdots,s(g^{n_\beta-1}(\beta))$ inside the triangulation. The non-degeneracy of $W$ was shown in \cite{Labardini-Fragoso2016}---in general, a potential $W$ for a quiver $Q$ is called non-degenerate if one can perform an arbitrary sequence of mutation of this quiver with potential, in the sense of \cite[Definition 7.2]{DerksenWeymanZelevinsky08}, without creating any 2-cycles.

\begin{lem}[{\cite[Prop 3.8, Lemma 4.1, Prop 4.2]{Ladkani2012-zx}}]\label{struc-jac-surf} The following hold in the completed Jacobian algebra $J$ of $(Q,W)$ with $Q=Q(T)$ and $W$ as in \eqref{eq:surface_potential}:
	\begin{enumerate}
		\item For any arrow $\alpha\in Q_1$, the elements $gf(\alpha)\cdot f(\alpha)\cdot \alpha$ and $fg(\alpha)\cdot g(\alpha)\cdot \alpha$ are equal to zero in $J$.
		\item The superfluous cycles $\alpha \cdot f^2(\alpha)\cdot f(\alpha)\cdot \alpha$ and $\alpha\cdot g^{n_\alpha-1}(\alpha)\cdots g(\alpha)\cdot \alpha$ are equal to zero in $J$. 
		\item The cycles $f^2(\alpha)\cdot f(\alpha)\cdot \alpha$ and $c_{\alpha}g^{n_\alpha-1}(\alpha)\cdots g(\alpha)\cdot \alpha$ (where $c_\alpha$ is the number associated to the puncture), as well as the expressions with $\alpha$ substituted by the other arrow $\overline{\alpha}$ with the same start $s(\overline{\alpha})=s(\alpha)$, are all equal in $J$. We denote their common value by $z_i$, where $i=s(\alpha)$.
		\item The Jacobian algebra $J$ is finite dimensional. It has a basis $$\{e_i\}_{i}\cup\{g^r(\alpha)\cdots g(\alpha)\cdot\alpha \}_{\alpha, r}\cup\{z_i\}_{i}$$ where $e_i$'s are the lazy paths, $\alpha\in Q_1, 0\le r\le n_{\alpha}-2$, and the $z_i$'s are defined above.
	\end{enumerate}
\end{lem}

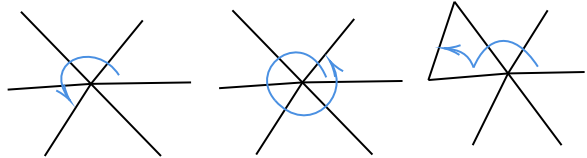
\begin{figure}[h]
    \centering

\tikzset{every picture/.style={line width=0.75pt}}         

\begin{tikzpicture}[x=0.75pt,y=0.75pt,yscale=-1,xscale=1,scale=0.75]

\draw [color={rgb, 255:red, 0; green, 0; blue, 0 }  ,draw opacity=1 ]   (85.03,86.79) -- (152.36,85.82) ;
\draw [color={rgb, 255:red, 0; green, 0; blue, 0 }  ,draw opacity=1 ]   (37.95,38.25) -- (85.03,86.79) ;
\draw [color={rgb, 255:red, 0; green, 0; blue, 0 }  ,draw opacity=1 ]   (85.03,86.79) -- (119.87,44.07) ;
\draw [color={rgb, 255:red, 0; green, 0; blue, 0 }  ,draw opacity=1 ]   (85.03,86.79) -- (53.96,135.33) ;
\draw [color={rgb, 255:red, 0; green, 0; blue, 0 }  ,draw opacity=1 ]   (85.03,86.79) -- (132.12,135.33) ;
\draw [color={rgb, 255:red, 0; green, 0; blue, 0 }  ,draw opacity=1 ]   (29,90.67) -- (85.03,86.79) ;
\draw [color={rgb, 255:red, 0; green, 0; blue, 0 }  ,draw opacity=1 ]   (227.23,85.82) -- (294.56,84.85) ;
\draw [color={rgb, 255:red, 0; green, 0; blue, 0 }  ,draw opacity=1 ]   (180.14,37.28) -- (227.23,85.82) ;
\draw [color={rgb, 255:red, 0; green, 0; blue, 0 }  ,draw opacity=1 ]   (227.23,85.82) -- (262.07,43.1) ;
\draw [color={rgb, 255:red, 0; green, 0; blue, 0 }  ,draw opacity=1 ]   (227.23,85.82) -- (196.15,134.36) ;
\draw [color={rgb, 255:red, 0; green, 0; blue, 0 }  ,draw opacity=1 ]   (227.23,85.82) -- (274.31,134.36) ;
\draw [color={rgb, 255:red, 0; green, 0; blue, 0 }  ,draw opacity=1 ]   (171.2,89.7) -- (227.23,85.82) ;
\draw [color={rgb, 255:red, 74; green, 144; blue, 226 }  ,draw opacity=1 ]   (103.87,80.96) .. controls (88.64,57.18) and (52.16,70.69) .. (69.76,96.37) ;
\draw [shift={(70.91,97.96)}, rotate = 232.83] [color={rgb, 255:red, 74; green, 144; blue, 226 }  ,draw opacity=1 ][line width=0.75]    (10.93,-3.29) .. controls (6.95,-1.4) and (3.31,-0.3) .. (0,0) .. controls (3.31,0.3) and (6.95,1.4) .. (10.93,3.29)   ;
\draw [color={rgb, 255:red, 74; green, 144; blue, 226 }  ,draw opacity=1 ]   (208.87,97.96) .. controls (189.09,75.14) and (228.64,46.98) .. (243.24,82.42) ;
\draw [color={rgb, 255:red, 74; green, 144; blue, 226 }  ,draw opacity=1 ]   (208.87,97.96) .. controls (230.98,123.03) and (260.85,95.69) .. (246.52,72.66) ;
\draw [shift={(245.59,71.26)}, rotate = 54.73] [color={rgb, 255:red, 74; green, 144; blue, 226 }  ,draw opacity=1 ][line width=0.75]    (10.93,-3.29) .. controls (6.95,-1.4) and (3.31,-0.3) .. (0,0) .. controls (3.31,0.3) and (6.95,1.4) .. (10.93,3.29)   ;
\draw [color={rgb, 255:red, 0; green, 0; blue, 0 }  ,draw opacity=1 ]   (365.5,79.77) -- (417,78.81) ;
\draw [color={rgb, 255:red, 0; green, 0; blue, 0 }  ,draw opacity=1 ]   (329.48,31.5) -- (365.5,79.77) ;
\draw [color={rgb, 255:red, 0; green, 0; blue, 0 }  ,draw opacity=1 ]   (365.5,79.77) -- (392.15,37.29) ;
\draw [color={rgb, 255:red, 0; green, 0; blue, 0 }  ,draw opacity=1 ]   (365.5,79.77) -- (341.73,128.04) ;
\draw [color={rgb, 255:red, 0; green, 0; blue, 0 }  ,draw opacity=1 ]   (365.5,79.77) -- (401.51,128.04) ;
\draw [color={rgb, 255:red, 0; green, 0; blue, 0 }  ,draw opacity=1 ]   (311.93,84.22) -- (365.5,79.77) ;
\draw [color={rgb, 255:red, 74; green, 144; blue, 226 }  ,draw opacity=1 ]   (342.33,75.81) .. controls (351.19,54.48) and (369.02,49.65) .. (385.56,75.23) ;
\draw [color={rgb, 255:red, 0; green, 0; blue, 0 }  ,draw opacity=1 ]   (311.93,84.22) -- (329.48,31.5) ;
\draw [color={rgb, 255:red, 74; green, 144; blue, 226 }  ,draw opacity=1 ]   (342.33,75.81) .. controls (338.74,65.91) and (332.16,64.04) .. (324.28,63.03) ;
\draw [shift={(322.34,62.8)}, rotate = 6.46] [color={rgb, 255:red, 74; green, 144; blue, 226 }  ,draw opacity=1 ][line width=0.75]    (10.93,-3.29) .. controls (6.95,-1.4) and (3.31,-0.3) .. (0,0) .. controls (3.31,0.3) and (6.95,1.4) .. (10.93,3.29)   ;

\end{tikzpicture}

    \caption{The leftmost is a nonzero path in $J$, while the right two equal zero.}
    \label{fig:illustration-paths}
\end{figure}

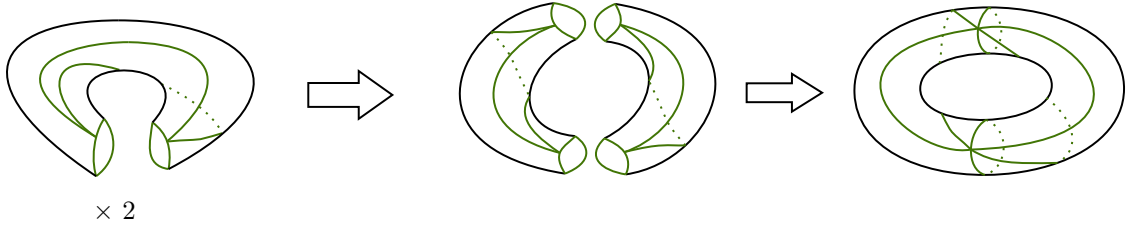
\begin{figure}[h]
    \centering

\tikzset{every picture/.style={line width=0.75pt}} 

\begin{tikzpicture}[x=0.75pt,y=0.75pt,yscale=-1,xscale=1]

\draw    (358.03,97.18) .. controls (294.8,87.15) and (296.05,47.81) .. (320.82,25.76) .. controls (329.19,18.3) and (340.27,12.82) .. (352.44,11.12) ;
\draw    (363.07,78.23) .. controls (347.06,76.67) and (339.86,68.19) .. (340.21,58.3) .. controls (340.56,48.14) and (348.86,36.48) .. (363.73,29.28) ;
\draw [color={rgb, 255:red, 65; green, 117; blue, 5 }  ,draw opacity=1 ]   (358.03,97.18) .. controls (356.79,94.94) and (354.04,90.57) .. (355.45,86.72) .. controls (356.87,82.86) and (362.16,80.08) .. (363.07,78.23) ;
\draw  [color={rgb, 255:red, 65; green, 117; blue, 5 }  ,draw opacity=1 ] (363.07,78.23) .. controls (370.92,83.36) and (372.04,91.65) .. (358.03,97.18) ;
\draw [color={rgb, 255:red, 65; green, 117; blue, 5 }  ,draw opacity=1 ]   (352.44,11.12) .. controls (350.53,14.73) and (350.77,19.08) .. (353.2,22.37) .. controls (355.62,25.66) and (360.02,27.72) .. (363.73,29.28) ;
\draw  [color={rgb, 255:red, 65; green, 117; blue, 5 }  ,draw opacity=1 ] (352.44,11.12) .. controls (365.24,13.09) and (373.97,19.81) .. (363.73,29.28) ;
\draw [color={rgb, 255:red, 65; green, 117; blue, 5 }  ,draw opacity=1 ]   (355.45,86.72) .. controls (300.51,72.51) and (323.09,31.85) .. (353.2,22.37) ;
\draw [color={rgb, 255:red, 65; green, 117; blue, 5 }  ,draw opacity=1 ]   (355.45,86.72) .. controls (346.42,80.4) and (332.12,70.53) .. (340.21,58.3) ;
\draw [color={rgb, 255:red, 65; green, 117; blue, 5 }  ,draw opacity=1 ]   (320.82,25.76) .. controls (330.62,25.53) and (333.63,27.11) .. (353.2,22.37) ;
\draw [color={rgb, 255:red, 65; green, 117; blue, 5 }  ,draw opacity=1 ] [dash pattern={on 0.84pt off 2.51pt}]  (320.82,25.76) .. controls (334.38,40.53) and (328.36,35.79) .. (340.21,58.3) ;
\draw    (383.38,11.52) .. controls (443.61,21.54) and (442.43,60.89) .. (418.83,82.94) .. controls (410.85,90.4) and (400.3,95.88) .. (388.7,97.57) ;
\draw    (378.58,30.46) .. controls (393.83,32.02) and (400.69,40.5) .. (400.36,50.39) .. controls (400.02,60.56) and (392.11,72.21) .. (377.94,79.41) ;
\draw [color={rgb, 255:red, 65; green, 117; blue, 5 }  ,draw opacity=1 ]   (383.38,11.52) .. controls (384.56,13.75) and (387.18,18.12) .. (385.83,21.98) .. controls (384.48,25.83) and (379.44,28.61) .. (378.58,30.46) ;
\draw  [color={rgb, 255:red, 65; green, 117; blue, 5 }  ,draw opacity=1 ] (378.58,30.46) .. controls (371.1,25.33) and (370.03,17.04) .. (383.38,11.52) ;
\draw [color={rgb, 255:red, 65; green, 117; blue, 5 }  ,draw opacity=1 ]   (388.7,97.57) .. controls (390.52,93.97) and (390.3,89.61) .. (387.98,86.32) .. controls (385.67,83.03) and (381.48,80.97) .. (377.94,79.41) ;
\draw  [color={rgb, 255:red, 65; green, 117; blue, 5 }  ,draw opacity=1 ] (388.7,97.57) .. controls (376.51,95.6) and (368.19,88.89) .. (377.94,79.41) ;
\draw [color={rgb, 255:red, 65; green, 117; blue, 5 }  ,draw opacity=1 ]   (385.83,21.98) .. controls (438.18,36.19) and (416.67,76.85) .. (387.98,86.32) ;
\draw [color={rgb, 255:red, 65; green, 117; blue, 5 }  ,draw opacity=1 ]   (385.83,21.98) .. controls (394.44,28.29) and (408.06,38.16) .. (400.36,50.39) ;
\draw [color={rgb, 255:red, 65; green, 117; blue, 5 }  ,draw opacity=1 ]   (418.83,82.94) .. controls (409.5,83.16) and (406.63,81.58) .. (387.98,86.32) ;
\draw [color={rgb, 255:red, 65; green, 117; blue, 5 }  ,draw opacity=1 ] [dash pattern={on 0.84pt off 2.51pt}]  (418.83,82.94) .. controls (405.91,68.16) and (411.65,72.9) .. (400.36,50.39) ;
\draw   (640,54.51) .. controls (640.27,27.25) and (606.15,13.31) .. (572.03,13.28) .. controls (537.49,13.26) and (502.96,27.48) .. (504.05,56.57) .. controls (505.11,84.7) and (537.15,97.94) .. (569.85,97.83) .. controls (582.54,97.79) and (595.32,95.74) .. (606.44,91.76) .. controls (625.61,84.91) and (639.82,72.34) .. (640,54.51) -- cycle ;
\draw   (603.57,52.6) .. controls (603.57,63.51) and (587,69.55) .. (570.43,70.05) .. controls (562.28,70.3) and (554.13,69.22) .. (547.94,66.72) .. controls (541.42,64.09) and (537.08,59.89) .. (537.22,54.03) .. controls (537.35,48.6) and (541.2,44.43) .. (546.99,41.53) .. controls (553.54,38.25) and (562.56,36.59) .. (571.48,36.56) .. controls (587.7,36.5) and (603.57,41.83) .. (603.57,52.6) -- cycle ;
\draw [color={rgb, 255:red, 65; green, 117; blue, 5 }  ,draw opacity=1 ]   (570.43,70.05) .. controls (565.53,72.99) and (562.94,79.2) .. (562.79,85.04) .. controls (562.63,90.92) and (564.94,96.43) .. (569.85,97.83) ;
\draw [color={rgb, 255:red, 65; green, 117; blue, 5 }  ,draw opacity=1 ] [dash pattern={on 0.84pt off 2.51pt}]  (570.43,70.05) .. controls (582.9,75.01) and (581.81,94.08) .. (569.85,97.83) ;
\draw [color={rgb, 255:red, 65; green, 117; blue, 5 }  ,draw opacity=1 ] [dash pattern={on 0.84pt off 2.51pt}]  (572.03,13.28) .. controls (580.73,20.18) and (579.64,33.53) .. (571.48,36.56) ;
\draw [color={rgb, 255:red, 65; green, 117; blue, 5 }  ,draw opacity=1 ]   (572.03,13.28) .. controls (568.81,15.8) and (566.95,19.89) .. (566.41,23.98) .. controls (565.68,29.54) and (567.4,35.09) .. (571.48,36.56) ;
\draw [color={rgb, 255:red, 65; green, 117; blue, 5 }  ,draw opacity=1 ]   (566.41,23.98) .. controls (475.23,35.91) and (532.87,90.26) .. (562.79,85.04) ;
\draw [color={rgb, 255:red, 65; green, 117; blue, 5 }  ,draw opacity=1 ]   (566.41,23.98) .. controls (614.44,14.46) and (671,81.2) .. (562.79,85.04) ;
\draw [color={rgb, 255:red, 65; green, 117; blue, 5 }  ,draw opacity=1 ]   (562.79,85.04) .. controls (552.45,75.01) and (554.08,79.3) .. (547.94,66.72) ;
\draw [color={rgb, 255:red, 65; green, 117; blue, 5 }  ,draw opacity=1 ]   (553.54,14.46) -- (566.41,23.98) ;
\draw [color={rgb, 255:red, 65; green, 117; blue, 5 }  ,draw opacity=1 ]   (566.41,23.98) -- (587.25,38.3) ;
\draw [color={rgb, 255:red, 65; green, 117; blue, 5 }  ,draw opacity=1 ]   (562.79,85.04) .. controls (574.2,92.17) and (584.53,92.17) .. (606.44,91.76) ;
\draw [color={rgb, 255:red, 65; green, 117; blue, 5 }  ,draw opacity=1 ] [dash pattern={on 0.84pt off 2.51pt}]  (546.99,41.53) .. controls (548.64,32.14) and (547.01,24.04) .. (553.54,14.46) ;
\draw [color={rgb, 255:red, 65; green, 117; blue, 5 }  ,draw opacity=1 ] [dash pattern={on 0.84pt off 2.51pt}]  (600.85,58.96) .. controls (616.62,71.24) and (617.71,85.54) .. (606.44,91.76) ;
\draw   (449.67,51.59) -- (472.51,51.59) -- (472.51,46.82) -- (487.74,56.35) -- (472.51,65.89) -- (472.51,61.12) -- (449.67,61.12) -- cycle ;
\draw [color={rgb, 255:red, 65; green, 117; blue, 5 }  ,draw opacity=1 ]   (121.59,98.49) .. controls (119.68,91.92) and (120.09,84.67) .. (121.61,78.78) .. controls (122.59,74.94) and (124.05,71.68) .. (125.65,69.53) ;
\draw  [color={rgb, 255:red, 65; green, 117; blue, 5 }  ,draw opacity=1 ] (125.65,69.53) .. controls (131.97,77.38) and (132.87,90.05) .. (121.59,98.49) ;
\draw [color={rgb, 255:red, 65; green, 117; blue, 5 }  ,draw opacity=1 ]   (150.01,71.34) .. controls (157.23,76.77) and (159.94,84.01) .. (158.13,94.87) ;
\draw  [color={rgb, 255:red, 65; green, 117; blue, 5 }  ,draw opacity=1 ] (150.01,71.34) .. controls (148.66,81) and (145.5,93.06) .. (158.13,94.87) ;
\draw [color={rgb, 255:red, 65; green, 117; blue, 5 }  ,draw opacity=1 ]   (133.32,44.8) .. controls (90.92,32.13) and (99,63.83) .. (121.61,78.78) ;
\draw [color={rgb, 255:red, 65; green, 117; blue, 5 }  ,draw opacity=1 ]   (157.23,81) .. controls (176.63,79.19) and (168.96,81.6) .. (185.2,76.77) ;
\draw [color={rgb, 255:red, 65; green, 117; blue, 5 }  ,draw opacity=1 ] [dash pattern={on 0.84pt off 2.51pt}]  (154.52,51.98) .. controls (167.6,58.01) and (180.69,67.12) .. (185.2,76.77) ;
\draw   (228.59,51.44) -- (254.03,51.44) -- (254.03,45.4) -- (271,57.47) -- (254.03,69.53) -- (254.03,63.5) -- (228.59,63.5) -- cycle ;
\draw    (125.65,69.53) .. controls (98.67,44.77) and (143.93,38.66) .. (154.52,51.98) .. controls (157.99,56.34) and (157.74,62.79) .. (150.01,71.34) ;
\draw    (121.59,98.49) .. controls (52,51.83) and (69,19.83) .. (133,19.83) ;
\draw    (133,19.83) .. controls (199,19.83) and (234,54.83) .. (158.13,94.87) ;
\draw [color={rgb, 255:red, 65; green, 117; blue, 5 }  ,draw opacity=1 ]   (121.61,78.78) .. controls (63,43.83) and (109,30.83) .. (138,30.83) ;
\draw [color={rgb, 255:red, 65; green, 117; blue, 5 }  ,draw opacity=1 ]   (157.23,81) .. controls (186,56.83) and (189,31.83) .. (138,30.83) ;

\draw (118.6,110.2) node [anchor=north west][inner sep=0.75pt]    {$\times \ 2$};

\end{tikzpicture}

    \caption{A surgery for constructing a $2:1$ covering} 
    \label{fig:surgery-covering}
\end{figure}

	Next, as in \cite[\S B]{MandelQin2021}, we consider coverings of a once-punctured surface $\?{\Sigma}=(\?{S},p)$ with triangulation $\?{T}$, and the associated covering quiver $Q(\?{T})$. Denote by $\pi:S\to \overline{S}$ a $d$-folded covering space map. The preimage of $\overline{T}$ in $S$ is a triangulation $T$ of $\Sigma=(S,\pi^{-1}(p))$.  Let $\bQ=Q(\?{T})$ and $Q=Q(T)$.   The induced covering map $\pi:\Sigma\rar \?{\Sigma}$ naturally induces a $d$-folded covering map of quivers $\pi:Q\to \bQ$ by sending the vertices and arrows of $Q$ to their counterparts in $\bQ$. See Figure \ref{fig:surgery-covering} for the construction of a $2$-folded covering space.

\begin{figure}[h]
    \centering
 \begin{tikzpicture}[xscale=0.6,yscale=0.75]
	\draw [color={rgb, 255:red, 65; green, 117; blue, 5 }] (0,0) rectangle (6,4);
	\draw [color={rgb, 255:red, 65; green, 117; blue, 5 } ]   (0,4) -- (6,0) ;
	\fill [color=blue](3,0)node[below]{$a$};\fill(3,4)[color=blue]node[above]{$a$};\fill[color=blue](0,2)node[left]{$b$};\fill[color=blue](6,1.2)node[left]{$e$};
	\fill [color=blue] (4,1.2)node[left]{$c$};
	\draw [<-] (3.5,2) -- (5.5,2);	\draw [<-] (3.15,3.5) -- (3.15,2.25);	\draw [<-]  (5.5,2.5) -- (3.5,3.5);
	\draw [<-] (2.5,2) -- (0.5,2); \draw[<-](2.85,0.5) -- (2.85,1.75);\draw [<-] (0.5,1.5)--(2.5,0.5);

    \draw [color={rgb, 255:red, 65; green, 117; blue, 5 }] (6,0) rectangle (12,4);
	\draw [color={rgb, 255:red, 65; green, 117; blue, 5 } ]   (6,4) -- (12,0) ;
	\fill [color=blue](9,0)node[below]{$d$};\fill(9,4)[color=blue]node[above]{$d$};\fill[color=blue](12,2)node[right]{$b$};
	\fill [color=blue] (10,1.2)node[left]{$f$};
	\draw [<-] (9.5,2) -- (11.5,2);	\draw [<-] (9.15,3.5) -- (9.15,2.25);	\draw [<-]  (11.5,2.5) -- (9.5,3.5);
	\draw [<-] (8.5,2) -- (6.5,2); \draw[<-](8.85,0.5) -- (8.85,1.75);\draw [<-] (6.5,1.5)--(8.5,0.5);
\end{tikzpicture}
    \caption{The quiver for the covering torus.
    }
    \label{fig:quiv-covering}
\end{figure}
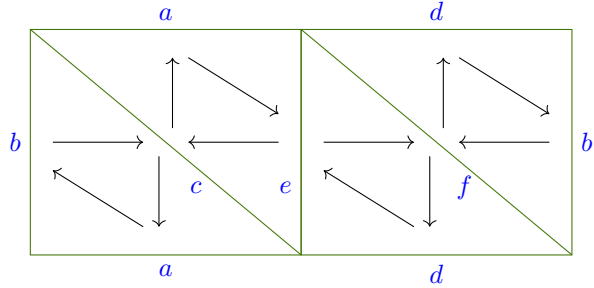

 We choose potentials $W$ for $Q$ and $\bW$ for $\bQ$ as in \eqref{eq:surface_potential} such that the complex numbers $c_p$ and $c_{\pi(p)}$ are the same for any puncture $p\in \Sigma$.

\begin{lem}\label{lem:compatible-potential-cov}
    If $\pi:\Sigma\to \?{\Sigma}$ is $d$-folded, then $\pi(W)=d\cdot \bW$ and $\sigma(\bW)=W$.
\end{lem}
\begin{proof}
    This follows from the properties of the $d$-folded covering map. 
    More precisely, the potentials consist of two parts: the enclosed triangles and rotations around punctures. For each triangle in the triangulation $\overline{T}$, its preimage consists of $d$-many triangles in $T$. Similarly, for any rotation of arcs in $T$ around a puncture $p\in \?{\Sigma}$, its preimage consists of the rotations around the $d$ punctures in $\pi^{-1}(p)$.
\end{proof}

\begin{figure}
    \centering
    
\begin{tikzpicture}
    [node distance=48pt,on grid,>={Stealth[round]},bend angle=45,   arc/.style={blue}, pre/.style={<-,shorten <=1pt,>={Stealth[round]},semithick},    post/.style={->,shorten >=1pt,>={Stealth[round]},semithick},  unfrozen/.style= {circle,inner sep=1pt,minimum size=12pt,draw=black!100},  frozen/.style={rectangle,inner sep=1pt,minimum size=12pt,draw=black!75,fill=cyan!100},   point/.style= {circle,inner sep=1pt,minimum size=5pt,draw=black!100,fill=black!100},   boundary/.style={-,draw=cyan},   internal/.style={-},    every label/.style= {black}]

     \node[unfrozen](p00) at (0,0){$m^{(0)}$};
    \node[unfrozen](p01) at (0,120pt){$m^{(0)}$};
    \node(b0) at (0,60pt){};
    \node(a00) at (60pt,0){};
    \node(a01) at (60pt,120pt){};
    \node(c0) at (60pt,60pt){};

     \node[unfrozen](p10) at (120pt,0){$m^{(1)}$};
    \node[unfrozen](p11) at (120pt,120pt){$m^{(1)}$};
    \node(b1) at (120pt,60pt){};
    \node(a10) at (180pt,0){};
    \node(a11) at (180pt,120pt){};
    \node(c1) at (180pt,60pt){};

     \node[unfrozen](p20) at (240pt,0){$m^{(2)}$};
    \node[unfrozen](p21) at (240pt,120pt){$m^{(2)}$};
    \node(b2) at (240pt,60pt){};
    \node(a20) at (300pt,0){};
    \node(a21) at (300pt,120pt){};
    \node(c2) at (300pt,60pt){};

     \node[unfrozen](p30) at (360pt,0){$m^{(0)}$};
    \node[unfrozen](p31) at (360pt,120pt){$m^{(0)}$};
    \node(b3) at (360pt,60pt){};

    \draw[arc] (p01) -- node[left, near end]{$b^{(0)}$} (p00);
    \draw[arc] (p11) -- node[left, near end]{$b^{(1)}$} (p10);
    \draw[arc] (p21) -- node[left, near end]{$b^{(2)}$} (p20);
    \draw[arc] (p31) -- node[right, near end]{$b^{(0)}$} (p30);

    \draw[arc] (p00) -- node [below]{$a^{(0)}$}(p10);
    \draw[arc] (p01) -- node [above]{$a^{(0)}$}(p11);
    \draw[arc] (p10) -- node [below]{$a^{(1)}$}(p20);
    \draw[arc] (p11) -- node [above]{$a^{(1)}$}(p21);
    \draw[arc] (p20) -- node [below]{$a^{(2)}$}(p30);
    \draw[arc] (p21) -- node [above]{$a^{(2)}$}(p31);

    \draw[arc] (p01) -- node [left, near start]{$c^{(0)}$} (p10);
    \draw[arc] (p11) -- node [left, near start]{$c^{(1)}$} (p20);
    \draw[arc] (p21) -- node [left, near start]{$c^{(2)}$} (p30);

    \draw [->][dashed] (a00) --node [left]{$\alpha_1^{(0)}$} (b0);
    \draw [->][dashed] (b0) --node [below]{$\beta_1^{(0)}$} (c0);
    \draw [->][dashed] (c0) --node [left]{$\gamma_1^{(0)}$} (a00);
    
    \draw [->][dashed] (a01) --node [right]{$\alpha_2^{(0)}$} (b1);
    \draw [->][dashed] (b1) --node [above]{$\beta_2^{(0)}$} (c0);
    \draw [->][dashed] (c0) --node [right]{$\gamma_2^{(0)}$} (a01);

    \draw [->][dashed] (a10) --node [left]{$\alpha_1^{(1)}$} (b1);
    \draw [->][dashed] (b1) --node [below]{$\beta_1^{(1)}$} (c1);
    \draw [->][dashed] (c1) --node [left]{$\gamma_1^{(1)}$} (a10);
    
    \draw [->][dashed] (a11) --node [right]{$\alpha_2^{(1)}$} (b2);
    \draw [->][dashed] (b2) --node [above]{$\beta_2^{(1)}$} (c1);
    \draw [->][dashed] (c1) --node [right]{$\gamma_2^{(1)}$} (a11);
    
    \draw [->][dashed] (a20) --node [left]{$\alpha_1^{(2)}$} (b2);
    \draw [->][dashed] (b2) --node [below]{$\beta_1^{(2)}$} (c2);
    \draw [->][dashed] (c2) --node [left]{$\gamma_1^{(2)}$} (a20);
    
    \draw [->][dashed] (a21) --node [right]{$\alpha_2^{(2)}$} (b3);
    \draw [->][dashed] (b3) --node [above]{$\beta_2^{(2)}$} (c2);
    \draw [->][dashed] (c2) --node [right]{$\gamma_2^{(2)}$} (a21);
    
    \end{tikzpicture}

    \caption{A triangulation of a $3:1$ covering $\Sigma$ of a once-punctured torus $\?{\Sigma}$. The nodes represent punctures, while solid arcs correspond to the vertices of the quiver. Nodes and arcs with the same labels are identified. The paths in the quiver are lifted to dashed paths on the surface. If we consider $\bQ$ as embedded in the middle part of the above quiver,  the labels of solid arcs also denote their gradings, with the exception of $b^{(0)}$ which does not live in $\supp_{Q_0}(P_k)$ at all.}
    \label{fig:triangulation_3_cover}
\end{figure}
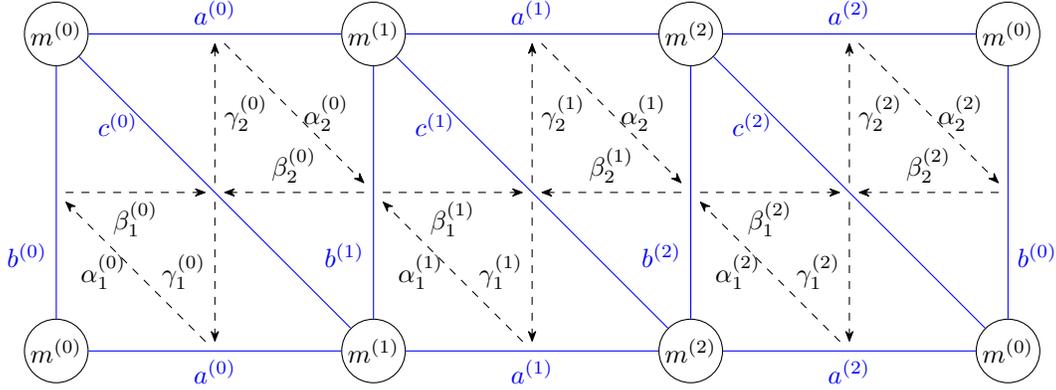

\begin{example}[Nice grading for a triple cover of the once-punctured torus]\label{eg:nice_grading_torus}
    Take $\?{\Sigma}$ to be a once-punctured torus. Let $\Sigma$ be a connected $3:1$ covering of $\?{\Sigma}$. This $\Sigma$ with the lifted triangulation is depicted in Figure \ref{fig:triangulation_3_cover}. Note that paths in $Q$ lift to dashed paths on $\Sigma$. The covering map $\pi$ sends $a^{(s)}$ to $a$ and similarly for $b^{(s)}$, $c^{(s)}$, where $s$ denotes the sheet label.
    
    We introduce the grading on $Q_0$ such that $\partial(i^{(s)})\coloneqq s$ for any arc $i^{(s)}$. Take $k=a^{(1)}$. We claim that $\partial$ restricts to a nice grading on $\supp_{Q_0}(P_k)$ with respect to $P_k$.  Note that the first condition of a nice grading is clearly satisfied.
    
    By Lemma \ref{struc-jac-surf}, the projective module $P_k$ has a basis consisting of the equivalence classes $[p]$, where $p$ takes the form:
    \begin{itemize}
        \item $e_{k}$
        \item $g^{r}(\alpha_1^{(1)})\cdots g(\alpha_1^{(1)})\cdot \alpha_1^{(1)}$ for $0\leq r\leq 5$, subpaths of an oriented dashed cycle around the puncture $m^{(1)}$
        \item $g^{r}(\alpha_2^{(1)})\cdots g(\alpha_2^{(1)})\cdot \alpha_2^{(1)}$ for $0\leq r\leq 5$, subpaths of an oriented dashed cycle around the puncture $m^{(2)}$
    \end{itemize}
    Therefore, $\supp_{Q_0}(P_k)=Q_0\setminus \{ b^{(0)}\}$.

    By Lemma \ref{struc-jac-surf}(1)(2), the only paths $p$ which represent nonzero elements in $P_k$ are the basis elements above, along with $\beta_1^{(1)}\alpha_1^{(1)}$, $\gamma_1^{(1)}\beta_1^{(1)}\alpha_1^{(1)}$, $\beta_2^{(1)}\alpha_2^{(1)}$, and $\gamma_2^{(1)}\beta_2^{(1)}\alpha_2^{(1)}$.  Thus,
    \begin{align*}
            \supp_{Q_1}(P_k)=\{\alpha_1^{(1)},\beta_1^{(1)},\gamma_1^{(1)},\alpha_2^{(1)},\beta_2^{(1)},\gamma_2^{(1)},\alpha_2^{(0)},\alpha_1^{(2)}, \beta_2^{(0)},\gamma_1^{(0)},\beta_1^{(2)},\gamma_2^{(2)}\}. 
    \end{align*}
    Introduce $\partial(\zeta)\coloneqq \partial(t(\zeta))-\partial(s(\zeta))$ for arrows $\zeta$.  It is straightforward to check that any pair of arrows in $\supp_{Q_1}(P_k)$ in the same $\Pi$-orbit will have the same degree---the key is that none of the arrows in $\supp_{Q_1}(P_k)$ go between $\partial^{-1}(0)$ and $\partial^{-1}(2)$.  The second condition for $\partial$ to be a nice grading follows.

    Note that while $\partial$ is $\Pi$-invariant on $\supp_{Q_1}(P_k)$, it is not a $\Pi$-invariant grading on $Q_1$, because $\partial(\alpha_2^{(0)})\neq\partial(\alpha_2^{(2)})$ and $\partial(\beta_2^{(0)})\neq\partial(\beta_2^{(2)})$. 
\end{example}

Example \ref{eg:nice_grading_torus} is generalized to prove the following.

\begin{lem}\label{lem:surface_nice_grading}
    Let $\Sigma$ be a connected $d:1$-covering of $\?{\Sigma}$ with $d\geq 3$ and deck transformation group $\Z/d\Z$.  Then for the associated $d:1$ covering of quivers with potential, there exists a nice grading with respect to $P_k$ for each $k$.
\end{lem}
\begin{proof}
    Recall that the vertices of $Q$ are arcs in $\Sigma$. The covering space $\Sigma$ is a union of $d$-sheets, where each sheet is obtained from $\?{\Sigma}$ by cutting an arc, which we denote by $\gamma$ (Figure \ref{fig:surgery-covering}). On $\Sigma$, we can label these sheets by $\Sigma^{(0)},\ldots,\Sigma^{(d-1)}$ consecutively. Each sheet $\Sigma^{(r)}$ contains two boundary arcs in $\pi^{-1}(\gamma)$, which will be denoted by $\gamma^{(r),-}$ and $\gamma^{(r),+}$ such that the latter is also contained in $\Sigma^{(r+1)}$ (denote $\Sigma^{(d)}=\Sigma^{(0)}$ here). Any vertex $i$ of $Q$ is an arc contained in some sheet $\Sigma^{(r)}$ (if $\pi(i)=\gamma$, we choose $\Sigma^{(r)}$ such that $i=\gamma^{(r),-}$). Define the grading $\partial(i)=r$. Without loss of generality, we can assume $\partial(k)=1$ and $\pi(k)\neq \gamma$. 
    
    Recall that, given any path (or any arrow) $p$ in $Q$, we can naturally lift $p$ to a curve in $\Sigma$, which intersects the relative interiors of the arcs consecutively as $p$ does. By abuses of notation, the curve is still called a path (or arrow) and denoted by $p$. Recall that $P_k$ has a basis consisting of equivalence classes $u_p$ of certain paths $p$ starting from $k$. We choose representatives $p$ as in Lemma \ref{struc-jac-surf}(4). Then $p$ appearing are subpaths of cycles surrounding an end point of the arc $k$. Moreover, it follows from Lemma \ref{struc-jac-surf}(1)(2) that the only other paths representing nonzero elements of $P_k$ are those of the form $f(a)\cdot a$ and $f^2(a)\cdot f(a)\cdot a$ for arrows $a$ starting from $k$.  
    We deduce that if the equivalence class of a path $q$ is nonzero and contained in $P_k$, then $q$ must be contained in $\Sigma^{(0)}\cup \Sigma^{(1)} \cup \Sigma^{(2)}\setminus (\gamma^{(0),-} \cup \gamma^{(2),+})$.

    Now, take any arrow $a\in \supp_{Q_1}(P_k)$. Then there is a path $q$ which contains $a$ and represents a nonzero element of $P_k$.  
    It follows that the arrow $a$ is contained in $\Sigma^{(0)}\cup \Sigma^{(1)} \cup \Sigma^{(2)}\setminus (\gamma^{(0),-} \cup \gamma^{(2),+})$. We define $\partial(a)=\partial(t(a))-\partial(s(a))$.

    Finally, let $g$ denote the deck transformation sending arcs $c$ in $\Sigma^{(r)}$ to $c'$ in $\Sigma^{(r+1)}$ such that $\pi(c')=\pi(c)$ (we denote $\Sigma^{(d)}=\Sigma^{(0)}$ here). This $g$ generates $\Pi$, and we see that the above grading on $\supp_{Q_1}(P_k)$ is invariant under the action of $g$, because there are no arrows $a\in \supp_{Q_1}(P_k)$ connecting arcs in $\Sigma^{(d-1)}$ to those in $\Sigma^{(0)}$. It follows that $\partial$ on $\supp_{Q_0}(P_k)$ is a nice grading.  
\end{proof}

\subsection{Proofs of Theorem \ref{thm:1p_surface_diagrams} and Theorem \ref{thm:theta-stability}}

\begin{proof}[Proof of Theorem \ref{thm:1p_surface_diagrams}]

Let $\Sigma'$ be the once-punctured surface obtained by reversing the orientation of $\Sigma$.  The triangulation $T$ of $\Sigma$ also gives a triangulation $T'$ of $\Sigma'$.  Let $Q(T')$ be the associated quiver and $W'$ the associated potential as in \eqref{eq:surface_potential}.  Then the associated quiver with potential is $J^{\op}$.  

Let $\wt{\Sigma'}$ be a connected $d:1$ covering of $\Sigma'$ with deck group $\Z/d\Z$, $d\geq 3$, equipped with the potential as in \eqref{eq:surface_potential}.  Let $\wt{J}^{\op}$ be the associated Jacobian algebra.  Then Lemma \ref{lem:surface_nice_grading} combined with Theorem \ref{thm:stability_diagram} tells us that $\?{\f{D}^{\st}(\wt{J}^{\op})}$ is equivalent to $\f{D}^{\st}(J^{\op})$.

Let $\wt{\Sigma}\rar \Sigma$ be the $d:1$ cover of $\Sigma$ obtained by reversing the orientation of $\wt{\Sigma'}$.  So $\wt{J}^{\op}$ is the opposite Jacobian algebra of $(\wt{\Sigma},\wt{W})$ for $\wt{W}$ as in \eqref{eq:surface_potential} for $\wt{\Sigma}$.  Since $\wt{\Sigma'}$ has more than one puncture, it is injective-reachable \cite{FominShapiroThurston08}, so \cite[Thm. 1.2.4]{qin2019bases} implies that $\frD^{\st}(\wt{J}^{\op})$ is equivalent to the cluster scattering diagram $\frD(\wt{\sd})$.

Finally, by \cite[\S B]{MandelQin2021}, the restriction $\?{\frD(\wt{\sd})}$ is equivalent to $\frD(\sd)$ if and only if $\Sigma$ is not a once-punctured torus. The desired claim follows.
\end{proof}

    \begin{rem}\label{rem:equivalence_coefficinets}
The proof of Theorem \ref{thm:1p_surface_diagrams} is also effective for principal coefficients. Therefore, $\f{D}^{\st,\prin}(J^{\op})$ and $\frD(\sd^{\prin})$ are equivalent if and only if $S$ is not a torus. It follows that the equivalence holds for arbitrary coefficients (i.e. an arbitrary choice of $B_{ij}$ such that not both of $i$ and $j$ are unfrozen). 
    \end{rem}

Note that the above proof of Theorem \ref{thm:1p_surface_diagrams} also shows the following:
\begin{lem}\label{lem:cover-stability}
    Let $\Sigma$ be a once-punctured closed triangulated surface.  Fix a potential $W$ as in \eqref{eq:surface_potential} for the associated quiver, and let $J^{\op}$ be the associated opposite Jacobian algebra.  Let $\pi:\wt{\Sigma}\rar \Sigma$ be a connected $d:1$ covering with deck group $\Z/d\Z$ and $d\geq 3$, and let $\wt{\sd}$ be the associated seed.  Then $\?{\frD(\wt{\sd})}$ is equivalent to $\frD^{\st}(J^{\op})$.
\end{lem}

\begin{proof}[Proof of Theorem \ref{thm:theta-stability}]
It was shown in \cite{MandelQin2021} that the bracelets basis and theta basis for once-punctured surfaces agree except in the case of the once-punctured torus if one uses the cluster scattering diagram for the surface.   Furthermore, it was shown that the bracelets basis and theta basis agree in all once-punctured closed surface cases if one uses $\?{\f{D}(\wt{\sd})}$ as in Lemma \ref{lem:cover-stability} in place of $\f{D}(\sd)$.  This fact plus Lemma \ref{lem:cover-stability} imply Theorem \ref{thm:theta-stability}.
\end{proof}

\section{Cyclic coverings with non-wrapping potentials}\label{sec:gen}

In this section, we apply ideas from the proof of Lemma \ref{lem:surface_nice_grading} to prove that the conclusion of Theorem \ref{thm:stability_diagram} holds under a more natural assumption than the existence of nice gradings.

Consider a $d:1$ covering $\pi:Q\rar \bQ$ with potentials $\bW=\pi(W_0)$, $W=\sigma(\bW)$, as usual.  Assume we have cyclic deck group $\Pi\cong \Z/d\Z$.  We can fix a decomposition $Q_0=\bigsqcup_{s=0}^{d-1} \bQ_0^{(s)}$ where each $\bQ_0^{(s)}$ is a copy of $\bQ_0$ and such that the action of any $k\in \Pi$ maps each $\bQ_0^{(i)}$ to $\bQ_0^{(i+k)}$ (superscripts taken mod $d$).  I.e., we view the sets $\bQ_0^{(i)}$ as the distinct sheets of vertices in our covering $\pi$.  For each $V\in Q_0$, define $\partial(V)=s$ where $V\in \bQ_0^{(s)}$.  Then for $\ba\in \bQ_1$ an arrow from $i$ to $j$, the lifts to $Q_1$ are arrows for each $s=0,1,\ldots,d-1$ from $i^{(s)}$ to $j^{(s+\delta_{\ba})}$ for some fixed $\delta_{\ba}\in \{0,1,\ldots,d-1\}$ (mod $d$).

In Example \ref{eg:nice_grading_torus} and the proof of Lemma \ref{lem:surface_nice_grading}, it was implicitly important that no cycle of $W$ wrapped all the way around the covering $\Sigma$.  We generalize this notion here.  Precisely, we say that $W$ is non-wrapping if we can find a function $\partial:\bQ_1\rightarrow \Z$ with $\partial(\ba)\in \{\delta_{\ba},\delta_{\ba}-d\}$ for each $\ba$ (think of this as choosing between viewing $a$ as going up $\delta_{\ba}$ sheets or going down $d-\delta_{\ba}$ sheets) such that, for each cycle $a_q\cdots a_1 a_0$ forming a term of $W$, we have $\sum_{i=0}^q \partial(\pi(a_i))=0$.  Intuitively, one can imagine collapsing $Q$ down to an oriented $d$-gon with each sheet corresponding to a vertex. 
 Choosing $\partial(\ba)$ amounts to viewing $a$ as going $\delta_{\ba}$ steps counterclockwise or $d-\delta_{\ba}$ steps clockwise around the $d$-gon.  Then $\frac{1}{d}\sum_{i=0}^q \partial(\pi(a_i))$ corresponds to the winding number of the cycle $a_q\cdots a_1a_0$. 

\begin{thm}\label{thm:non-wrap}
    Let $\pi:(Q,W)\rar (\bQ,\bW)$ be a covering of quivers with potential with deck group $\Pi=\Z/d\Z$ and sheets labelled compatibly, where $d> 1$.  Suppose $W$ is non-wrapping.  Then $\?{\f{D}^{\st}(J)}$ is equivalent to $\f{D}^{\st}(\bJ)$.
\end{thm}

\begin{proof}
For any $l>1$, we can extend $\pi$ to a $2ld:1$ covering $\wt{\pi}:\wt{Q}\rightarrow \?{Q}$ with potential $\wt{W}=\wt{\sigma}(\bW)$ for $\wt{\sigma}$ associated to $\wt{\pi}$. The vertices are $\wt{Q}_0=\bigsqcup_{s=0}^{2l d-1} \bQ_0^{(s)}$, and the arrows are determined by specifying that the preimage of an arrow $\ba$ from $i$ to $j$ consists of, for each $s=0,...,2l d-1$, an arrow from $i^{(s)}$ to $j^{(s+\partial(\ba))}$ for $(s+\partial(\ba))$ interpreted modulo $2ld$.  The condition that $W$ is non-wrapping ensures that the terms of $W$ (hence of $\bW$) really do lift to cycles in $\wt{Q}$, so the lift $\wt{W}=\wt{\sigma}(\bW)$ genuinely gives a potential for $\wt{Q}$.

Let $k\in \wt{Q}_0$.  Without loss of generality, we may assume $k\in \?{Q}_0^{(ld)}$ (otherwise we relabel the sheets).  Let $\wt{P}_k$ denote the associated projective module. We claim that the map $\partial:\wt{Q}_0\rightarrow \{0,1,\ldots,2ld-1\}\subset \ZZ$, $\partial(\bQ_0^{(s)})=s$, restricts to a nice grading on $\supp_{\wt{Q}_0}(\wt{P}_k^l)$ with respect to $\wt{P}_k^l$ (cf. Remark \ref{rmk:nice-Pl}).  Recall, this means that:
\begin{enumerate}
    \item for each $i,i'\in \supp_{\wt{Q}_0}(\wt{P}_k^l)$, we have $\partial(i)=\partial(i')$ if and only if $i=i'$, and
    \item for any arrows $a,a'\in \supp_{\wt{Q}_1}(\wt{P}_k^l)$ with $\pi(a)=\pi(a')$, we have $\partial(t(a))-\partial(s(a))=\partial(t(a'))-\partial(s(a'))$.
\end{enumerate}  
Condition (1) is clearly satisfied.  Furthermore, since $0\leq |\partial(\ba)|\leq d-1$, and since  $l$ is at least the length of the longest nonzero path in $\wt{P}_k^l$, every arrow $a$ in $\supp_{\wt{Q}_1}(\wt{P}_k^l)$ satisfies $\partial(t(a))-\partial(s(a))=\partial(\wt{\pi}(a))$, so Condition (2) holds as well.

As indicated in Remark \ref{rmk:nice-Pl}, the proof of Corollary \ref{cor:projection_Euler} applies to show that $\chi(\Quot_{\bn}(\bP_{\bk}))=\sum_{\wt{n}\in \wt{\pi}^{-1}(\bn)} \chi(\Quot_{\wt{n}}(\wt{P}_{\wt{k}}))$ whenever $|n|\leq l+1$.  It then follows as in the proof of Theorem \ref{thm:stability_diagram} that $\?{\f{D}^{\st}(\wt{J})}$ and $\f{D}^{\st}(\bJ)$ are equivalent to order $l$, i.e., $\?{\f{D}^{\st}(\wt{J})}_{l}$ and $\f{D}^{\st}(\bJ)_{l}$ are equivalent.

Furthermore, we note that $\wt{\pi}:\wt{Q}\rar \bQ$ factors through the original $\pi:Q\rar \bQ$, giving a covering map $\pi':\wt{Q}\rar Q$ with $\sigma'(W)=\wt{W}$ for $\sigma'$ associated to $\pi'$.  A similar argument shows that $\?{\f{D}^{\st}(\wt{J})}'_{l}$ (the folded scattering diagram with respect to $\pi'$) is equivalent to $\f{D}^{\st}(J)_{l}$.  The restriction construction in \S \ref{sec:restriction} respects compositions of coverings, so we find that in fact $\?{\f{D}^{\st}(J)}_{l}$ is equivalent to $\f{D}^{\st}(\bJ)_{l}$.  Finally, since $l$ was arbitrary, we deduce that $\?{\f{D}^{\st}(J)}$ is equivalent to $\f{D}^{\st}(\bJ)$.
\end{proof}

\begin{rem}\label{rem:solvable}
    Theorem \ref{thm:non-wrap} applies to cyclic coverings, but one may apply the result iteratively to obtain equivalences of scattering diagrams for more general (e.g., solvable) deck groups as well.
\end{rem}

\begin{eg}\label{eg:general-manifold-covering}
    Let $(\bQ,\bW)$ be a quiver with potential.  Suppose we have an immersion of (the topological realization of) $\bQ$ into a surface $\?{\SSS}$ (or in fact, any manifold with boundary $\?{\SSS}$) such that, for each term of $\bW$, the associated loop in $\?{\SSS}$ is contractible.  Let $\pi:\SSS\rar \?{\SSS}$ be a $d:1$ normal covering space of $\?{\SSS}$ with cyclic deck group.  Then $\pi$ determines a $d:1$ covering $\pi:(Q,W)\rar (\bQ,\bW)$ of quivers with potential.  By choosing a branch cut for the covering $\pi:\SSS\rar \?{\SSS}$, we can identify $\SSS$ as a union of $d$ sheets, each containing a copy $\bQ^{(s)}_0$ of $\bQ_0$ (labelled sequentially and cyclically with $s=0,\ldots,d-1$).  The assumption that terms of $\bW$ correspond to contractible loops in $\?{\SSS}$ ensures that $W$ is non-wrapping.  We can thus apply Theorem \ref{thm:non-wrap} to conclude that $\?{\f{D}^{\st}(J)}$ is equivalent to $\f{D}^{\st}(\bJ)$.

    By applying this construction iteratively as in Remark \ref{rem:solvable}, this easily generalizes from cyclic deck groups to solvable deck groups.
\end{eg}

\begin{eg}\label{eg:Lie-Grass}
Consider the following quiver $\bQ$ (Figure \ref{fig:another_base_quiver}) which arises in the study of cluster algebras from Lie theory and affine Grassmanians (see \cite{GeissLeclercSchroeer10} \cite[Figure 6.1]{cautis2019cluster})


\begin{figure}[h]
    \centering

\tikzset{every picture/.style={line width=0.75pt}} 

\begin{tikzpicture}[x=0.75pt,y=0.75pt,yscale=-1,xscale=1]

\draw (144.05,163.45) node [anchor=north west][inner sep=0.75pt]    {$A$};
\draw (394.94,164.67) node [anchor=north west][inner sep=0.75pt]    {$C$};
\draw (268.64,59.69) node [anchor=north west][inner sep=0.75pt]    {$B$};
\draw (506.31,59.69) node [anchor=north west][inner sep=0.75pt]    {$D$};
\draw (185.23,94.42) node [anchor=north west][inner sep=0.75pt]    {$a_{1}$};
\draw (269.25,153.31) node [anchor=north west][inner sep=0.75pt]    {$c$};
\draw (393.11,70.61) node [anchor=north west][inner sep=0.75pt]    {$d$};
\draw (227.68,121.17) node [anchor=north west][inner sep=0.75pt]    {$a_{2}$};
\draw (349.58,101.87) node [anchor=north west][inner sep=0.75pt]    {$b_{1}$};
\draw (309.83,125.83) node [anchor=north west][inner sep=0.75pt]    {$b_{2}$};
\draw (431.21,98.87) node [anchor=north west][inner sep=0.75pt]    {$c_{1}$};
\draw (481.87,119.73) node [anchor=north west][inner sep=0.75pt]    {$c_{2}$};
\draw    (503.31,67.29) -- (285.64,67.29) ;
\draw [shift={(283.64,67.29)}, rotate = 360] [color={rgb, 255:red, 0; green, 0; blue, 0 }  ][line width=0.75]    (10.93,-3.29) .. controls (6.95,-1.4) and (3.31,-0.3) .. (0,0) .. controls (3.31,0.3) and (6.95,1.4) .. (10.93,3.29)   ;
\draw    (406.92,160.27) -- (501.86,71.19) ;
\draw [shift={(503.31,69.82)}, rotate = 136.82] [color={rgb, 255:red, 0; green, 0; blue, 0 }  ][line width=0.75]    (10.93,-3.29) .. controls (6.95,-1.4) and (3.31,-0.3) .. (0,0) .. controls (3.31,0.3) and (6.95,1.4) .. (10.93,3.29)   ;
\draw    (410.94,170.21) -- (506.37,80.66) ;
\draw [shift={(507.83,79.29)}, rotate = 136.82] [color={rgb, 255:red, 0; green, 0; blue, 0 }  ][line width=0.75]    (10.93,-3.29) .. controls (6.95,-1.4) and (3.31,-0.3) .. (0,0) .. controls (3.31,0.3) and (6.95,1.4) .. (10.93,3.29)   ;
\draw    (283.64,68.25) -- (393.24,158.99) ;
\draw [shift={(394.78,160.27)}, rotate = 219.62] [color={rgb, 255:red, 0; green, 0; blue, 0 }  ][line width=0.75]    (10.93,-3.29) .. controls (6.95,-1.4) and (3.31,-0.3) .. (0,0) .. controls (3.31,0.3) and (6.95,1.4) .. (10.93,3.29)   ;
\draw    (281.3,79.29) -- (390.4,169.62) ;
\draw [shift={(391.94,170.89)}, rotate = 219.62] [color={rgb, 255:red, 0; green, 0; blue, 0 }  ][line width=0.75]    (10.93,-3.29) .. controls (6.95,-1.4) and (3.31,-0.3) .. (0,0) .. controls (3.31,0.3) and (6.95,1.4) .. (10.93,3.29)   ;
\draw    (391.94,172.22) -- (163.05,171.11) ;
\draw [shift={(161.05,171.1)}, rotate = 0.28] [color={rgb, 255:red, 0; green, 0; blue, 0 }  ][line width=0.75]    (10.93,-3.29) .. controls (6.95,-1.4) and (3.31,-0.3) .. (0,0) .. controls (3.31,0.3) and (6.95,1.4) .. (10.93,3.29)   ;
\draw    (157.57,159.05) -- (264.11,69.61) ;
\draw [shift={(265.64,68.32)}, rotate = 139.99] [color={rgb, 255:red, 0; green, 0; blue, 0 }  ][line width=0.75]    (10.93,-3.29) .. controls (6.95,-1.4) and (3.31,-0.3) .. (0,0) .. controls (3.31,0.3) and (6.95,1.4) .. (10.93,3.29)   ;
\draw    (161.05,169.18) -- (266.59,80.58) ;
\draw [shift={(268.13,79.29)}, rotate = 139.99] [color={rgb, 255:red, 0; green, 0; blue, 0 }  ][line width=0.75]    (10.93,-3.29) .. controls (6.95,-1.4) and (3.31,-0.3) .. (0,0) .. controls (3.31,0.3) and (6.95,1.4) .. (10.93,3.29)   ;

\end{tikzpicture}

    \caption{The base quiver $\bQ$}
    \label{fig:another_base_quiver}
\end{figure}

Consider as well the potential $\bW=cb_1a_1+cb_2a_2+dc_1b_1+dc_2b_2$.  Then the Jacobian algebra $\bJ$ is the one considered in \cite[Theorem 6.4]{BuanIyamaReitenSmith08} where $\bJ$ is known to be finite-dimensional and $\bW$ is non-degenerate.

Consider the following $d$-layer quiver $Q$, $d\geq 2$ (Figure \ref{fig:another_covering_quiver}). Notice that we identify the top and bottom layer.  This $Q$ naturally admits a $d:1$ covering map $\pi:Q\to \bQ$---we take lifts of $a_1,b_1,c_1,c,d$ to be horizontal, while lifts of $a_2,b_2,c_2$ go up or down one sheet.  The induced potential $W=\sigma(\bW)$ is given by summing all the triangles having $\overrightarrow{C_iA_i}$ or $\overrightarrow{D_iB_i}$ as an edge, where $i$ ranges from $1$ to $d$. It is non-wrapping.


\begin{figure}[h]
    \centering

\tikzset{every picture/.style={line width=0.75pt}} 

\begin{tikzpicture}[x=0.75pt,y=0.75pt,yscale=-1,xscale=1]

\draw (142.92,68.54) node [anchor=north west][inner sep=0.75pt]    {$A_{1}$};
\draw (313.12,68.89) node [anchor=north west][inner sep=0.75pt]    {$C_{1}$};
\draw (227.68,39.01) node [anchor=north west][inner sep=0.75pt]    {$B_{1}$};
\draw (388.45,39.01) node [anchor=north west][inner sep=0.75pt]    {$D_{1}$};
\draw (142.92,121.02) node [anchor=north west][inner sep=0.75pt]    {$A_{2}$};
\draw (313.12,121.37) node [anchor=north west][inner sep=0.75pt]    {$C_{2}$};
\draw (227.68,91.49) node [anchor=north west][inner sep=0.75pt]    {$B_{2}$};
\draw (388.45,91.49) node [anchor=north west][inner sep=0.75pt]    {$D_{2}$};
\draw (136.92,170.07) node [anchor=north west][inner sep=0.75pt]    {$A_{3}$};
\draw (307.12,170.42) node [anchor=north west][inner sep=0.75pt]    {$C_{3}$};
\draw (221.68,140.54) node [anchor=north west][inner sep=0.75pt]    {$B_{3}$};
\draw (382.45,140.54) node [anchor=north west][inner sep=0.75pt]    {$D_{3}$};
\draw (135.92,332.66) node [anchor=north west][inner sep=0.75pt]    {$A_{1}$};
\draw (306.12,333.03) node [anchor=north west][inner sep=0.75pt]    {$C_{1}$};
\draw (220.68,301.2) node [anchor=north west][inner sep=0.75pt]    {$B_{1}$};
\draw (381.45,301.2) node [anchor=north west][inner sep=0.75pt]    {$D_{1}$};
\draw (143.92,280.42) node [anchor=north west][inner sep=0.75pt]    {$A_{d}$};
\draw (314.12,280.79) node [anchor=north west][inner sep=0.75pt]    {$C_{d}$};
\draw (228.68,248.96) node [anchor=north west][inner sep=0.75pt]    {$B_{d}$};
\draw (389.45,248.96) node [anchor=north west][inner sep=0.75pt]    {$D_{d}$};
\draw    (385.45,47.61) -- (251.68,47.61) ;
\draw [shift={(249.68,47.61)}, rotate = 360] [color={rgb, 255:red, 0; green, 0; blue, 0 }  ][line width=0.75]    (10.93,-3.29) .. controls (6.95,-1.4) and (3.31,-0.3) .. (0,0) .. controls (3.31,0.3) and (6.95,1.4) .. (10.93,3.29)   ;
\draw    (336.12,72.4) -- (383.58,53.82) ;
\draw [shift={(385.45,53.09)}, rotate = 158.62] [color={rgb, 255:red, 0; green, 0; blue, 0 }  ][line width=0.75]    (10.93,-3.29) .. controls (6.95,-1.4) and (3.31,-0.3) .. (0,0) .. controls (3.31,0.3) and (6.95,1.4) .. (10.93,3.29)   ;
\draw    (249.68,51.96) -- (308.23,72.31) ;
\draw [shift={(310.12,72.97)}, rotate = 199.17] [color={rgb, 255:red, 0; green, 0; blue, 0 }  ][line width=0.75]    (10.93,-3.29) .. controls (6.95,-1.4) and (3.31,-0.3) .. (0,0) .. controls (3.31,0.3) and (6.95,1.4) .. (10.93,3.29)   ;
\draw    (310.12,77.46) -- (169.92,77.18) ;
\draw [shift={(167.92,77.17)}, rotate = 0.12] [color={rgb, 255:red, 0; green, 0; blue, 0 }  ][line width=0.75]    (10.93,-3.29) .. controls (6.95,-1.4) and (3.31,-0.3) .. (0,0) .. controls (3.31,0.3) and (6.95,1.4) .. (10.93,3.29)   ;
\draw    (167.92,72.18) -- (222.8,52.71) ;
\draw [shift={(224.68,52.04)}, rotate = 160.47] [color={rgb, 255:red, 0; green, 0; blue, 0 }  ][line width=0.75]    (10.93,-3.29) .. controls (6.95,-1.4) and (3.31,-0.3) .. (0,0) .. controls (3.31,0.3) and (6.95,1.4) .. (10.93,3.29)   ;
\draw    (385.45,100.09) -- (251.68,100.09) ;
\draw [shift={(249.68,100.09)}, rotate = 360] [color={rgb, 255:red, 0; green, 0; blue, 0 }  ][line width=0.75]    (10.93,-3.29) .. controls (6.95,-1.4) and (3.31,-0.3) .. (0,0) .. controls (3.31,0.3) and (6.95,1.4) .. (10.93,3.29)   ;
\draw    (336.12,124.88) -- (383.58,106.3) ;
\draw [shift={(385.45,105.57)}, rotate = 158.62] [color={rgb, 255:red, 0; green, 0; blue, 0 }  ][line width=0.75]    (10.93,-3.29) .. controls (6.95,-1.4) and (3.31,-0.3) .. (0,0) .. controls (3.31,0.3) and (6.95,1.4) .. (10.93,3.29)   ;
\draw    (249.68,104.43) -- (308.23,124.79) ;
\draw [shift={(310.12,125.44)}, rotate = 199.17] [color={rgb, 255:red, 0; green, 0; blue, 0 }  ][line width=0.75]    (10.93,-3.29) .. controls (6.95,-1.4) and (3.31,-0.3) .. (0,0) .. controls (3.31,0.3) and (6.95,1.4) .. (10.93,3.29)   ;
\draw    (310.12,129.94) -- (169.92,129.65) ;
\draw [shift={(167.92,129.65)}, rotate = 0.12] [color={rgb, 255:red, 0; green, 0; blue, 0 }  ][line width=0.75]    (10.93,-3.29) .. controls (6.95,-1.4) and (3.31,-0.3) .. (0,0) .. controls (3.31,0.3) and (6.95,1.4) .. (10.93,3.29)   ;
\draw    (167.92,124.65) -- (222.8,105.19) ;
\draw [shift={(224.68,104.52)}, rotate = 160.47] [color={rgb, 255:red, 0; green, 0; blue, 0 }  ][line width=0.75]    (10.93,-3.29) .. controls (6.95,-1.4) and (3.31,-0.3) .. (0,0) .. controls (3.31,0.3) and (6.95,1.4) .. (10.93,3.29)   ;
\draw    (379.45,149.14) -- (245.68,149.14) ;
\draw [shift={(243.68,149.14)}, rotate = 360] [color={rgb, 255:red, 0; green, 0; blue, 0 }  ][line width=0.75]    (10.93,-3.29) .. controls (6.95,-1.4) and (3.31,-0.3) .. (0,0) .. controls (3.31,0.3) and (6.95,1.4) .. (10.93,3.29)   ;
\draw    (330.12,173.93) -- (377.58,155.35) ;
\draw [shift={(379.45,154.62)}, rotate = 158.62] [color={rgb, 255:red, 0; green, 0; blue, 0 }  ][line width=0.75]    (10.93,-3.29) .. controls (6.95,-1.4) and (3.31,-0.3) .. (0,0) .. controls (3.31,0.3) and (6.95,1.4) .. (10.93,3.29)   ;
\draw    (243.68,153.48) -- (302.23,173.84) ;
\draw [shift={(304.12,174.5)}, rotate = 199.17] [color={rgb, 255:red, 0; green, 0; blue, 0 }  ][line width=0.75]    (10.93,-3.29) .. controls (6.95,-1.4) and (3.31,-0.3) .. (0,0) .. controls (3.31,0.3) and (6.95,1.4) .. (10.93,3.29)   ;
\draw    (304.12,178.99) -- (163.92,178.7) ;
\draw [shift={(161.92,178.7)}, rotate = 0.12] [color={rgb, 255:red, 0; green, 0; blue, 0 }  ][line width=0.75]    (10.93,-3.29) .. controls (6.95,-1.4) and (3.31,-0.3) .. (0,0) .. controls (3.31,0.3) and (6.95,1.4) .. (10.93,3.29)   ;
\draw    (161.92,173.7) -- (216.8,154.24) ;
\draw [shift={(218.68,153.57)}, rotate = 160.47] [color={rgb, 255:red, 0; green, 0; blue, 0 }  ][line width=0.75]    (10.93,-3.29) .. controls (6.95,-1.4) and (3.31,-0.3) .. (0,0) .. controls (3.31,0.3) and (6.95,1.4) .. (10.93,3.29)   ;
\draw    (378.45,309.8) -- (244.68,309.8) ;
\draw [shift={(242.68,309.8)}, rotate = 360] [color={rgb, 255:red, 0; green, 0; blue, 0 }  ][line width=0.75]    (10.93,-3.29) .. controls (6.95,-1.4) and (3.31,-0.3) .. (0,0) .. controls (3.31,0.3) and (6.95,1.4) .. (10.93,3.29)   ;
\draw    (329.12,336.21) -- (376.6,316.41) ;
\draw [shift={(378.45,315.64)}, rotate = 157.37] [color={rgb, 255:red, 0; green, 0; blue, 0 }  ][line width=0.75]    (10.93,-3.29) .. controls (6.95,-1.4) and (3.31,-0.3) .. (0,0) .. controls (3.31,0.3) and (6.95,1.4) .. (10.93,3.29)   ;
\draw    (242.68,314.43) -- (301.24,336.12) ;
\draw [shift={(303.12,336.81)}, rotate = 200.32] [color={rgb, 255:red, 0; green, 0; blue, 0 }  ][line width=0.75]    (10.93,-3.29) .. controls (6.95,-1.4) and (3.31,-0.3) .. (0,0) .. controls (3.31,0.3) and (6.95,1.4) .. (10.93,3.29)   ;
\draw    (303.12,341.6) -- (162.92,341.29) ;
\draw [shift={(160.92,341.29)}, rotate = 0.13] [color={rgb, 255:red, 0; green, 0; blue, 0 }  ][line width=0.75]    (10.93,-3.29) .. controls (6.95,-1.4) and (3.31,-0.3) .. (0,0) .. controls (3.31,0.3) and (6.95,1.4) .. (10.93,3.29)   ;
\draw    (160.92,335.97) -- (215.81,315.23) ;
\draw [shift={(217.68,314.53)}, rotate = 159.31] [color={rgb, 255:red, 0; green, 0; blue, 0 }  ][line width=0.75]    (10.93,-3.29) .. controls (6.95,-1.4) and (3.31,-0.3) .. (0,0) .. controls (3.31,0.3) and (6.95,1.4) .. (10.93,3.29)   ;
\draw    (386.45,257.56) -- (252.68,257.56) ;
\draw [shift={(250.68,257.56)}, rotate = 360] [color={rgb, 255:red, 0; green, 0; blue, 0 }  ][line width=0.75]    (10.93,-3.29) .. controls (6.95,-1.4) and (3.31,-0.3) .. (0,0) .. controls (3.31,0.3) and (6.95,1.4) .. (10.93,3.29)   ;
\draw    (337.12,283.97) -- (384.6,264.17) ;
\draw [shift={(386.45,263.4)}, rotate = 157.37] [color={rgb, 255:red, 0; green, 0; blue, 0 }  ][line width=0.75]    (10.93,-3.29) .. controls (6.95,-1.4) and (3.31,-0.3) .. (0,0) .. controls (3.31,0.3) and (6.95,1.4) .. (10.93,3.29)   ;
\draw    (250.68,262.19) -- (309.24,283.88) ;
\draw [shift={(311.12,284.57)}, rotate = 200.32] [color={rgb, 255:red, 0; green, 0; blue, 0 }  ][line width=0.75]    (10.93,-3.29) .. controls (6.95,-1.4) and (3.31,-0.3) .. (0,0) .. controls (3.31,0.3) and (6.95,1.4) .. (10.93,3.29)   ;
\draw    (311.12,289.36) -- (170.92,289.05) ;
\draw [shift={(168.92,289.05)}, rotate = 0.13] [color={rgb, 255:red, 0; green, 0; blue, 0 }  ][line width=0.75]    (10.93,-3.29) .. controls (6.95,-1.4) and (3.31,-0.3) .. (0,0) .. controls (3.31,0.3) and (6.95,1.4) .. (10.93,3.29)   ;
\draw    (168.92,283.73) -- (223.81,262.99) ;
\draw [shift={(225.68,262.29)}, rotate = 159.31] [color={rgb, 255:red, 0; green, 0; blue, 0 }  ][line width=0.75]    (10.93,-3.29) .. controls (6.95,-1.4) and (3.31,-0.3) .. (0,0) .. controls (3.31,0.3) and (6.95,1.4) .. (10.93,3.29)   ;
\draw    (167.12,116.62) -- (223.26,61.32) ;
\draw [shift={(224.68,59.92)}, rotate = 135.44] [color={rgb, 255:red, 0; green, 0; blue, 0 }  ][line width=0.75]    (10.93,-3.29) .. controls (6.95,-1.4) and (3.31,-0.3) .. (0,0) .. controls (3.31,0.3) and (6.95,1.4) .. (10.93,3.29)   ;
\draw    (249.68,59.59) -- (308.68,116.12) ;
\draw [shift={(310.12,117.51)}, rotate = 223.78] [color={rgb, 255:red, 0; green, 0; blue, 0 }  ][line width=0.75]    (10.93,-3.29) .. controls (6.95,-1.4) and (3.31,-0.3) .. (0,0) .. controls (3.31,0.3) and (6.95,1.4) .. (10.93,3.29)   ;
\draw    (335.17,116.97) -- (386.04,62.08) ;
\draw [shift={(387.4,60.61)}, rotate = 132.82] [color={rgb, 255:red, 0; green, 0; blue, 0 }  ][line width=0.75]    (10.93,-3.29) .. controls (6.95,-1.4) and (3.31,-0.3) .. (0,0) .. controls (3.31,0.3) and (6.95,1.4) .. (10.93,3.29)   ;
\draw    (161.92,166.35) -- (223.18,112.41) ;
\draw [shift={(224.68,111.09)}, rotate = 138.64] [color={rgb, 255:red, 0; green, 0; blue, 0 }  ][line width=0.75]    (10.93,-3.29) .. controls (6.95,-1.4) and (3.31,-0.3) .. (0,0) .. controls (3.31,0.3) and (6.95,1.4) .. (10.93,3.29)   ;
\draw    (249.68,112.43) -- (302.7,164.78) ;
\draw [shift={(304.12,166.18)}, rotate = 224.64] [color={rgb, 255:red, 0; green, 0; blue, 0 }  ][line width=0.75]    (10.93,-3.29) .. controls (6.95,-1.4) and (3.31,-0.3) .. (0,0) .. controls (3.31,0.3) and (6.95,1.4) .. (10.93,3.29)   ;
\draw    (330.12,166.55) -- (384.44,114.47) ;
\draw [shift={(385.89,113.09)}, rotate = 136.21] [color={rgb, 255:red, 0; green, 0; blue, 0 }  ][line width=0.75]    (10.93,-3.29) .. controls (6.95,-1.4) and (3.31,-0.3) .. (0,0) .. controls (3.31,0.3) and (6.95,1.4) .. (10.93,3.29)   ;
\draw    (160.92,328.42) -- (224.21,270.38) ;
\draw [shift={(225.68,269.03)}, rotate = 137.48] [color={rgb, 255:red, 0; green, 0; blue, 0 }  ][line width=0.75]    (10.93,-3.29) .. controls (6.95,-1.4) and (3.31,-0.3) .. (0,0) .. controls (3.31,0.3) and (6.95,1.4) .. (10.93,3.29)   ;
\draw    (250.24,270.56) -- (302.71,327.16) ;
\draw [shift={(304.07,328.63)}, rotate = 227.17] [color={rgb, 255:red, 0; green, 0; blue, 0 }  ][line width=0.75]    (10.93,-3.29) .. controls (6.95,-1.4) and (3.31,-0.3) .. (0,0) .. controls (3.31,0.3) and (6.95,1.4) .. (10.93,3.29)   ;
\draw    (329.12,328.67) -- (385.99,271.98) ;
\draw [shift={(387.4,270.56)}, rotate = 135.09] [color={rgb, 255:red, 0; green, 0; blue, 0 }  ][line width=0.75]    (10.93,-3.29) .. controls (6.95,-1.4) and (3.31,-0.3) .. (0,0) .. controls (3.31,0.3) and (6.95,1.4) .. (10.93,3.29)   ;
\draw  [dash pattern={on 4.5pt off 4.5pt}]  (162,276.02) -- (223.14,163.89) ;
\draw [shift={(224.1,162.14)}, rotate = 118.6] [color={rgb, 255:red, 0; green, 0; blue, 0 }  ][line width=0.75]    (10.93,-3.29) .. controls (6.95,-1.4) and (3.31,-0.3) .. (0,0) .. controls (3.31,0.3) and (6.95,1.4) .. (10.93,3.29)   ;
\draw  [dash pattern={on 4.5pt off 4.5pt}]  (239.8,162.14) -- (314.4,274.72) ;
\draw [shift={(315.51,276.39)}, rotate = 236.47] [color={rgb, 255:red, 0; green, 0; blue, 0 }  ][line width=0.75]    (10.93,-3.29) .. controls (6.95,-1.4) and (3.31,-0.3) .. (0,0) .. controls (3.31,0.3) and (6.95,1.4) .. (10.93,3.29)   ;
\draw  [dash pattern={on 4.5pt off 4.5pt}]  (330.55,276.39) -- (386.13,163.93) ;
\draw [shift={(387.02,162.14)}, rotate = 116.3] [color={rgb, 255:red, 0; green, 0; blue, 0 }  ][line width=0.75]    (10.93,-3.29) .. controls (6.95,-1.4) and (3.31,-0.3) .. (0,0) .. controls (3.31,0.3) and (6.95,1.4) .. (10.93,3.29)   ;

\end{tikzpicture}

    \caption{The covering quiver $Q$}
    \label{fig:another_covering_quiver}
\end{figure}

It is straightforward to see that Theorem \ref{thm:non-wrap} applies to this setting.  Hence, for $J$ and $\bJ$ the Jacobian algebras associated to $(Q,W)$ and $(\bQ,\bW)$, respectively, we have that $\?{\f{D}^{\st}(J)}$ is equivalent to $\f{D}^{\st}(\bJ)$.

\end{eg}

		
        \def\cprime{$'$}
\providecommand{\bysame}{\leavevmode\hbox to3em{\hrulefill}\thinspace}
\providecommand{\MR}{\relax\ifhmode\unskip\space\fi MR }
\providecommand{\MRhref}[2]{%
  \href{http://www.ams.org/mathscinet-getitem?mr=#1}{#2}
}
\providecommand{\href}[2]{#2}

\end{document}